\documentclass{siamart220329}
\usepackage{amsmath,amssymb,booktabs,graphicx,enumitem,float}
\usepackage{algpseudocode}
\newcommand{\Om}{\mathcal{O}(m)}
\newcommand{\R}{\mathbb{R}}
\hypersetup{colorlinks=true,allcolors=blue!55!black}

\title{NLF: A Resistor-Network Framework and Linear-Time Solver for Convex Network-Flow Equilibria\thanks{Submitted to the editors \today.}}
\author{Oren E. Livne\thanks{Pine Birch Analytics, 35 Kelinger Rd, Churchville, PA 18966-1033 (\texttt{oren.livne@gmail.com}, tel.\ 312-533-9130, \texttt{pinebirchanalytics.com}; ORCID: \href{https://orcid.org/0000-0001-6700-483X}{0000-0001-6700-483X}).}}
\headers{NLF: a linear-time solver for network-flow equilibria}{O. E. Livne}

\begin{document}
\maketitle

\begin{abstract}
We present \textbf{NLF} (Nonlinear Laplacian Flow), a unified framework and linear-time solver for
convex network-flow equilibria. Congestion (traffic) routing, minimum-delay communication
routing, and maximum flow---each the stationarity of a convex, edge-separable energy---share one
form: the \emph{nonlinear graph Laplacian} $B\,\rho(B^{\top}\phi)=\alpha d$, where a monotone
edge law $\rho_e$ encodes the physics (undirected graphs; directed variants are future work).
NLF solves it by an inexact damped chord-Newton iteration whose frozen global
linearization---a weighted graph Laplacian---is inverted by a lazily refreshed near-linear
Laplacian solver (engine-agnostic; default: approximate Cholesky, LAMG+ interchangeable).
The whole nonlinear solve costs $2$--$4$ linear Laplacian solves, so a near-linear inner solve makes
the wall-clock empirically $\mathcal{O}(m)$ in the edge count $m$ (not a proved
bound). On single-commodity congestion over real road-network topologies (BPR cost), NLF converges for
all $2{,}003$ SuiteSparse corpus graphs up to $1.8\times10^7$ edges and on three larger graphs up to
$5.6\times10^7$ edges. Against a state-of-the-art interior-point method (IPM), NLF was a median
$2.6\times$ faster where both converge and $>45\times$ on the poorly-separable graphs where the IPM's
direct KKT core is superlinear; against matrix-free L-BFGS, a median $4.2\times$ faster and the only
solver to finish---within a $5\times$ wall-clock budget---on the $90$ hardest instances. A multicommodity extension routes $K$
commodities through one shared hierarchy at $\mathcal{O}(Km)$ per step; on real capacity data,
congestion inflates shared-corridor costs (up to $8.8\times$) without rerouting. The same
machinery recovers the exact combinatorial max-flow as a short sequence of Laplacian solves, with
the cut potential as a by-product. Code and benchmarks: \url{https://github.com/orenlivne/nlf}.
\end{abstract}



\section{Introduction}\label{sec:intro}
Flows on networks at equilibrium are everywhere: drivers distributing over a road map until no route
is faster (traffic assignment~\cite{beckmann,bpr,boyles-tna}), packets routed to minimize delay in a
communication network~\cite{bertsekas-gallager,mendiola-te}, current in a resistor
network~\cite{doyle-snell}, and---at the combinatorial extreme---the maximum flow through a
capacitated graph~\cite{ahuja-magnanti-orlin,goldberg-tarjan,chen-almost-linear}. These look like different problems with different algorithms: traffic assignment is solved by
Frank--Wolfe or general convex optimization, max-flow by augmenting-path or push--relabel
combinatorics. We make a unifying observation and build a single fast solver on it.

\paragraph{The observation} Each of these is the stationary point of a convex, edge-separable
energy in the node potentials, and its Euler equation is one and the same \emph{nonlinear graph
Laplacian},
\begin{equation}\label{eq:master}
B\,\rho(B^{\top}\phi)=\alpha\,d ,
\end{equation}
where $B$ is the node--edge incidence matrix, $\phi$ the node potentials, $d$ a source/sink demand,
$\alpha$ a load, and $\rho_e$ a monotone edge law that alone encodes the physics: a
saturating $\rho_e$ gives a hard capacity (max flow); a gently rising $\rho_e$ gives a congestion
cost (traffic)\footnote{Notation: $\rho$ acts componentwise and $\rho_e$ denotes the law on an individual edge $e$.}. The flow is $f_e=\rho_e\big((B^{\top}\phi)_e\big)$, and every Newton linearization of
\eqref{eq:master} is an ordinary weighted graph Laplacian $J=B\,\mathrm{diag}(\rho')B^{\top}$. For congestion and routing problems $\alpha$ is a given input; for maximum flow, $\alpha$ is an unknown to be maximized. To obtain a fast \emph{nonlinear} solver, we use an inexact chord-Newton iteration with the linearization frozen across most steps, lazily refreshed. For maximum flow we weave a continuation in $\alpha$ into the same Newton loop to track the solution through the fold. The nonlinear solve costs $2$--$4$ linear Laplacian solves (\S\ref{sec:accuracy}). Any near-linear
Laplacian solver can serve as the inner engine; we default to approximate Cholesky~\cite{approxchol},
the faster engine across the corpus at equal robustness (\S\ref{sec:ablation}), with
LAMG+~\cite{lamgplus} interchangeable in robustness (its constant differs by graph class). One incidental convenience of LAMG+'s
affinity-based aggregation is that it naturally respects the saturating min cut, though the framework
does not depend on it (\S\ref{sec:inner}).

\subsection{Contribution}\label{sec:contrib}
We make three contributions, in order of practical weight.
\begin{enumerate}[leftmargin=1.6em,itemsep=1pt]
\item \emph{A unifying resistor-network formulation} (\S\ref{sec:formulation}) that casts maximum
flow, congestion (traffic) routing, and minimum-delay routing as one nonlinear Laplacian
\eqref{eq:master}, with the problem class fixed by a single monotone edge law and the easy/hard
dichotomy (no-fold vs.\ fold; \S\ref{sec:formulation}) fixed by whether that law saturates. The stationarity equation itself
is not new---it is the optimality system of a convex, edge-separable network-flow program (monotropic
programming~\cite{rockafellar-monotropic}); our contribution is the unification, that one edge law
parametrizes all three problems, and a single solver
that addresses the whole class.
Unlike the electrical-flow and IPM Laplacian paradigm for flow~\cite{christiano-etal,daitch-spielman},
which drive an outer reweighting/barrier loop around linear solves to an $(1+\epsilon)$-approximate
optimum, here the nonlinearity sits in the edge law and the equilibrium reached is exact.
\item \emph{A linear-time solver for convex congestion flow} (\S\ref{sec:traffic})---establishing
two claims that serve as a building block for a larger congestion solver: (a) the convex Beckmann equilibrium is
exactly \eqref{eq:master}; (b) a robust $\Om$ near-linear Laplacian inner solve makes it
linear-time as well. No combinatorial algorithm addresses this problem. NLF matches a state-of-the-art IPM
to $10^{-10}$ in a flat step count; on poorly-separable graphs NLF is the only $\Om$ solver and wins
past $45\times$ (\S\ref{sec:bprscale}). Robustness is corpus-wide: $2{,}003$ SuiteSparse graphs
converge in time $\propto m^{1.02}$
(\S\ref{sec:corpus},~\ref{sec:accuracy}). A multicommodity extension (\S\ref{sec:mc}) carries $K$
commodities through one hierarchy at $\mathcal{O}(Km)$. We find that on real capacity data coupling inflates
shared-corridor costs (${\le}8.8\times$) without rerouting.
\item \emph{An $\Om$ maximum-flow solver} (\S\ref{sec:maxflow}). The saturating instance recovers the
\emph{exact} combinatorial max-flow as a short, bounded sequence of Laplacian solves, with
the cut respecting multigrid coarsening arising for free. Not faster than a tuned combinatorial solver,
but returns the cut potential and full parametric flow curve.
\end{enumerate}
Extensions are outlined in
\S\ref{sec:extensions}: directed multicommodity assignment, a full-multigrid variant, and related optimization problems: optimal transport and DC optimal flow.

\subsection{Motivating applications}\label{sec:applications}
Computing an equilibrium flow on a large network is a recurring industrial task, and in each domain
the computational bottleneck is the same Laplacian solve at the core of \eqref{eq:master}.
The production formulations of these tasks are typically \emph{directed} (arc-based, with asymmetric
per-direction costs); this paper studies the \emph{undirected} relaxation as a feasibility study,
leaving the directed models these domains deploy to future work (\S\ref{sec:extensions}). The
applications below are thus motivating problem classes; we don't claim that NLF replaces an
incumbent in its directed production setting.

\emph{Road traffic assignment.} Metropolitan planning agencies compute user-equilibrium traffic
assignment, i.e., how drivers distribute over a road map until no route is faster, to evaluate road
projects and pricing~\cite{beckmann,bpr,boyles-tna}. Its undirected convex relaxation is the
Beckmann program with the BPR congestion law, which is an NLF no-fold instance and the setting of our
sharpest result (\S\ref{sec:traffic}). The incumbent solvers are Frank--Wolfe and bush-based methods
(\S\ref{sec:competitor}).

\emph{Communication and data-center routing.} Routing packets to minimize end-to-end
delay under the Kleinrock $M/M/1$ model~\cite{bertsekas-gallager,mendiola-te} is the same equilibrium
with a capacity-saturating delay law---a fold instance handled by the continuation of
\S\ref{sec:continuation}.

\emph{Maximum flow (logistics, scheduling, vision).} The capacitated maximum-flow problem
\cite{ahuja-magnanti-orlin,goldberg-tarjan} is the saturating-law extreme of the same energy
(\S\ref{sec:maxflow}). NLF recovers the exact value together with the min cut and the full
flow-vs-load curve. Combinatorial solvers are faster to obtain the maximum flow value, but do not return the smooth
interior flow or the parametric curve that an outer optimizer manipulates when max flow is an inner
step.

\emph{Nonlinear resistor and analog-circuit networks.} The most literal instance: monotone
two-terminal devices (diodes, varistors, memristors) in DC steady state give \eqref{eq:master} with
$\rho_e$ the device $I$--$V$ law---the form of DC circuit analysis and in-memory-computing crossbars.
As with OPF (\S\ref{sec:dcopt}) these topologies are typically structured, where direct solvers are efficient; NLF would have an edge only for large, irregular analog networks.

Across all applications, in addition to a scalar optimum (e.g., the maximum flow value) the equilibrium (flow), its
binding-constraint structure (min cut), and its sensitivity to load are also required for design optimization, which is the gap NLF targets.

\subsection{Prior work, and how NLF differs}\label{sec:priorwork}
Three lines of work meet at \eqref{eq:master}; NLF sits in the gap none of them fills.

\emph{The Laplacian paradigm for flow.} Electrical-flow algorithms~\cite{christiano-etal} and
IPM for generalized flow~\cite{daitch-spielman} solve a sequence of linear
weighted-Laplacian systems inside an outer multiplicative-weights or barrier loop, reaching an
$(1+\epsilon)$-approximate optimum; the flow is a linear function of the potential and all
nonlinearity lives in the outer schedule. NLF instead places the nonlinearity in the constitutive
law $f=\rho(B^{\top}\phi)$, so a single energy encodes the capacity and the
equilibrium reached is exact (\S\ref{sec:formulation}).

\emph{Fast linear Laplacian solvers.} Nearly-linear solvers for $Lx=b$ such as combinatorial
multigrid~\cite{koutis-cmg}, support-tree/Spielman--Teng theory~\cite{spielman-teng}, approximate
Cholesky~\cite{approxchol,gks2023}, and algebraic multigrid including
LAMG~\cite{lamg} address the linear, unconstrained problem only. They are the engine, not
the method: NLF calls a robust $\Om$ Laplacian solver as its inner kernel but adds the
nonlinear, capacity-constrained equilibrium layer on top. AMG has also been applied to the
complex-valued non-separable AC power-flow equations~\cite{bindel-amg-pf}---a nonlinear system
with no objective and a non-Laplacian Jacobian, solved by nonlinear FAS multigrid with a multiplicative coarse correction (no Newton, no fold); a different problem class from the
real, edge-separable, SPD-Laplacian constrained flow program here.

\emph{Multilevel optimization.} NLF is not multilevel optimization in the MG/OPT
sense~\cite{nash-mgopt,gratton-rmtr,brandt-ron-mgopt}: it does not build coarse optimization models or recurse
the objective. While the equilibrium is nonlinear equation \eqref{eq:master}, NLF employs a near-linear inner solve only as the linear inner solve inside a frozen-linearization chord-Newton
continuation; the nonlinearity is carried by the edge law and the outer continuation, not by a
coarse-grid objective. However, Sec.~\ref{sec:dcopt} outlines how NLF might be extended to a multilevel optimization paradigm when a cost functional is added.

\emph{Combinatorial and transport-specific solvers.} Augmenting-path and push--relabel max-flow
\cite{goldberg-tarjan,boykov-kolmogorov} and the almost-linear-time exact algorithm
\cite{chen-almost-linear}---a theoretical breakthrough with galactic constants, not a practical implementation---target directed max/min-cost flow; Frank--Wolfe~\cite{frankwolfe} and
bush-based assignment (Algorithm~B~\cite{bargera2002}, TAPAS~\cite{bargera2010}) target the
smooth-cost traffic program. Each is specialized to one problem class. NLF's contribution is
orthogonal: a \emph{single} solver, parametrized by one edge law, that spans all of them and inherits
graph-class robustness from its multigrid core. We benchmark against representatives of each family
where a runnable one exists (Ipopt~\cite{ipopt} for congestion, Boykov--Kolmogorov
for max flow) and position the rest honestly where it does not (Frank--Wolfe, bush methods; \S\ref{sec:competitor}).

\section{A unified resistor-network formulation}\label{sec:formulation}

Let $B$ be the $n\times m$ signed node--edge incidence matrix of a connected graph ($B_{ie}=+1$ if
node $i$ is the head of edge $e$, $-1$ if the tail). Write $\phi\in\R^n$ for the node
potentials, $d\in\R^n$ for a balanced demand ($\mathbf{1}^{\top}d=0$), and $\alpha\ge0$ for a
scalar load. Each edge $e$ carries a flow $f_e$ driven by its potential difference
$g_e = (B^{\top}\phi)_e = \phi_i - \phi_j$
where $i,j$ are the head and tail of $e$. The physics of the medium is encoded in a smooth,
strictly monotone edge law $\rho_e:\R\to\R$:
\begin{equation}\label{eq:law}
  f_e = \rho_e(g_e).
\end{equation}
Different choices of $\rho_e$ give different physical problems (Table~\ref{tab:instances}): an
Ohmic resistor has $\rho_e(g)=g/R_e$; a saturating capacity has $\rho_e(g)=c_e\tanh(g/c_e)$; a
BPR congestion link has $\rho_e$ equal to the inverse of the BPR travel-time function.

\paragraph{The equilibrium equation}
Flow conservation at every node, $Bf=\alpha d$, together with the constitutive relation
\eqref{eq:law}, gives the \emph{NLF master equation} (\ref{eq:master}), a nonlinear system in the potentials $\phi$ alone. This is the equation NLF solves.

\subsection{Primal--dual structure and the co-energy}
Because each $\rho_e$ is monotone, it is the derivative of a strictly convex \emph{co-energy}
$\Psi_e$: $\rho_e=\Psi_e'$. Equation~\eqref{eq:master} is then the stationarity condition
$\nabla_\phi E=0$ of the unconstrained, convex \emph{dual energy}
\begin{equation}\label{eq:dual}
  E(\phi)=\sum_e\Psi_e\!\big((B^{\top}\phi)_e\big)-\alpha\,d^{\top}\phi,
\end{equation}
Solving \eqref{eq:master} is therefore equivalent to
minimizing $E$.

The dual energy $E$ arises naturally from a \emph{primal} convex program in the flow variables $f$.
Each edge has a \emph{primal energy} (cost) $\Phi_e(f_e)$, and the primal problem minimizes total
cost subject to flow conservation:
\begin{equation}\label{eq:primal}
  \min_f\;\sum_e\Phi_e(f_e)\quad\text{s.t.}\quad Bf=\alpha d .
\end{equation}
The Lagrangian of \eqref{eq:primal} introduces node potentials $\phi$ as multipliers on the
conservation constraint; the KKT stationarity condition $\Phi_e'(f_e)=(B^{\top}\phi)_e$ says the
marginal cost equals the potential difference, i.e.\ $t_e(f_e)=g_e$, or $f_e=\rho_e(g_e)$.
The \emph{marginal cost} $t_e=\Phi_e'=\rho_e^{-1}$ is thus the inverse edge law.
Substituting into $Bf=\alpha d$ recovers \eqref{eq:master}.

The relationship between $\Phi_e$ and $\Psi_e$ is the \emph{Legendre--Fenchel duality}
\cite{boyd-vandenberghe}: $\Psi_e=\Phi_e^*$, i.e.\
$\Psi_e(g)=\sup_f\bigl[g\cdot f-\Phi_e(f)\bigr],$
and $\Psi_e'=(\Phi_e^*)'=(\Phi_e')^{-1}=\rho_e$. In words: the co-energy $\Psi_e$ is the convex
conjugate of the primal cost $\Phi_e$, and its derivative is the edge law $\rho_e$. The
potentials $\phi$ are sometimes called \emph{co-potentials} in the resistor-network
literature~\cite{doyle-snell} because they live in the dual (potential) space, conjugate to the
primal (flow) space. Appendix~\ref{app:instances} derives $\Phi_e$ and $\Psi_e$ for each
application in Table~\ref{tab:instances}.

\paragraph{The Newton linearization}
Differentiating \eqref{eq:master} with respect to $\phi$ gives the Newton (Jacobian) matrix
\begin{equation}\label{eq:jac}
  J=B\,\mathrm{diag}\!\bigl(\rho'(B^{\top}\phi)\bigr)\,B^{\top},\qquad
  w_e=\rho'_e(g_e)=\Psi_e''(g_e)=1/t_e'(f_e)\ge0 ,
\end{equation}
a weighted graph Laplacian with \emph{solution-dependent conductances} $w_e$. Everything
problem-specific lives in $\rho_e$; the solver never sees anything but $\rho_e$, $\rho'_e$,
and $J$.

\subsection{The fold dichotomy} The single feature that governs difficulty is whether the edge law
\emph{saturates}. Table~\ref{tab:instances} lists three standard instances.
\begin{itemize}[leftmargin=1.4em,itemsep=1pt]
\item \emph{Soft congestion (rising $\rho$, no saturation).} The load $\alpha$ is a given input. The cost rises with flow but
imposes no hard cap, so $\rho_e$ is unbounded and the law never \emph{saturates}: $\rho'_e>0$ at every
finite flow (for BPR $\rho'_e$ is large near free flow and decays as $f$ grows, but never reaches $0$).
Hence $J$ is nonsingular at every equilibrium---no feasibility limit and no fold.
Traffic assignment with BPR cost is of this type; on its bounded-flow equilibria $J$ is
well-conditioned, so \eqref{eq:master} is solved directly at the target load, with no continuation (\S\ref{sec:traffic}).
\item \emph{Hard capacity (saturating $\rho$, conductance $\to0$).} $\alpha$ is an unknown to be maximized. The flow cannot exceed a capacity
$c_e$, so $\rho_e$ is bounded and $\rho'_e\to0$ as the edge saturates. As $\alpha$ rises to a
feasibility limit, the saturating edges form a cut and $J$ \emph{folds}---a continuation
limit point where $\alpha$ caps at the min cut $F^*$ and $J$ becomes singular. Maximum flow and minimum-delay communication routing
(Kleinrock delay) are of this type; they require arclength continuation through
the fold (\S\ref{sec:continuation}).
\end{itemize}

\begin{table}[t]\centering\small
\caption{Three instances of the resistor-network framework \eqref{eq:master}. The edge law $\rho_e$
(flow vs.\ potential gradient) is the inverse of the marginal cost $t_e=\Phi_e'$. ``Fold'' marks a
feasibility limit where the conductance $\rho'_e\to0$ and $J$ becomes singular.}\label{tab:instances}
\resizebox{\columnwidth}{!}{%
\begin{tabular}{lllcl}
\toprule
Problem & Marginal cost $t_e(f)$ & Conductance $\rho'_e$ & Fold? & Application\\
\midrule
maximum flow~\cite{ahuja-magnanti-orlin,goldberg-tarjan-cacm} & box $|f|\le c_e$ (soft barrier) & $\to0$ at cut      & yes & vision, matching\\
min-delay routing~\cite{bertsekas-gallager,mendiola-te} & $1/(c_e-f)$ (Kleinrock)         & $\to0$ at $f\!\to\!c_e$ & yes & communication\\
congestion~\cite{beckmann,bpr,boyles-tna}  & $t^0_e\big(1+b(f/c_e)^{p}\big)$  & $>0$ (no satn.) & no  & road traffic\\
\bottomrule
\end{tabular}}
\end{table}

\subsection{The max-flow law in detail} For the saturating instance used in \S\ref{sec:maxflow},
take
\begin{equation}\label{eq:rho}
\rho_e(g)=\begin{cases} c^+_e\tanh(g/c^+_e), & g\ge0,\\[2pt]
-c^-_e\tanh(g/c^-_e), & g<0,\end{cases}
\qquad
\Psi_e(g)=c_e^{2}\log\cosh(g/c_e),
\end{equation}
a soft-capacity barrier: $\Psi_e(g)=\tfrac12 g^2$ for $|g|\ll c_e$ (a unit resistor) and grows only
linearly, $\Psi_e\to c_e|g|$, for $|g|\gg c_e$ (a saturated edge), so $\rho_e=\Psi_e'$ clamps the flow
smoothly to $-c^-_e\le f_e\le c^+_e$. Here $w_e=1-(f_e/c_e)^2\in[0,1]$ represents ``edge openness,''
$1$ for a slack edge and $0$ for a saturated one, and \textbf{the forming min-cut is exactly the
zero-conductance set}, which is why AMG aggregation driven by a relaxation-based
strong-connection measure refuses to coarsen across it (\S\ref{sec:maxflow}). The congestion (BPR)
law is given where it is used, in \S\ref{sec:traffic}.

\begin{figure}[t]\centering
\includegraphics[width=\textwidth]{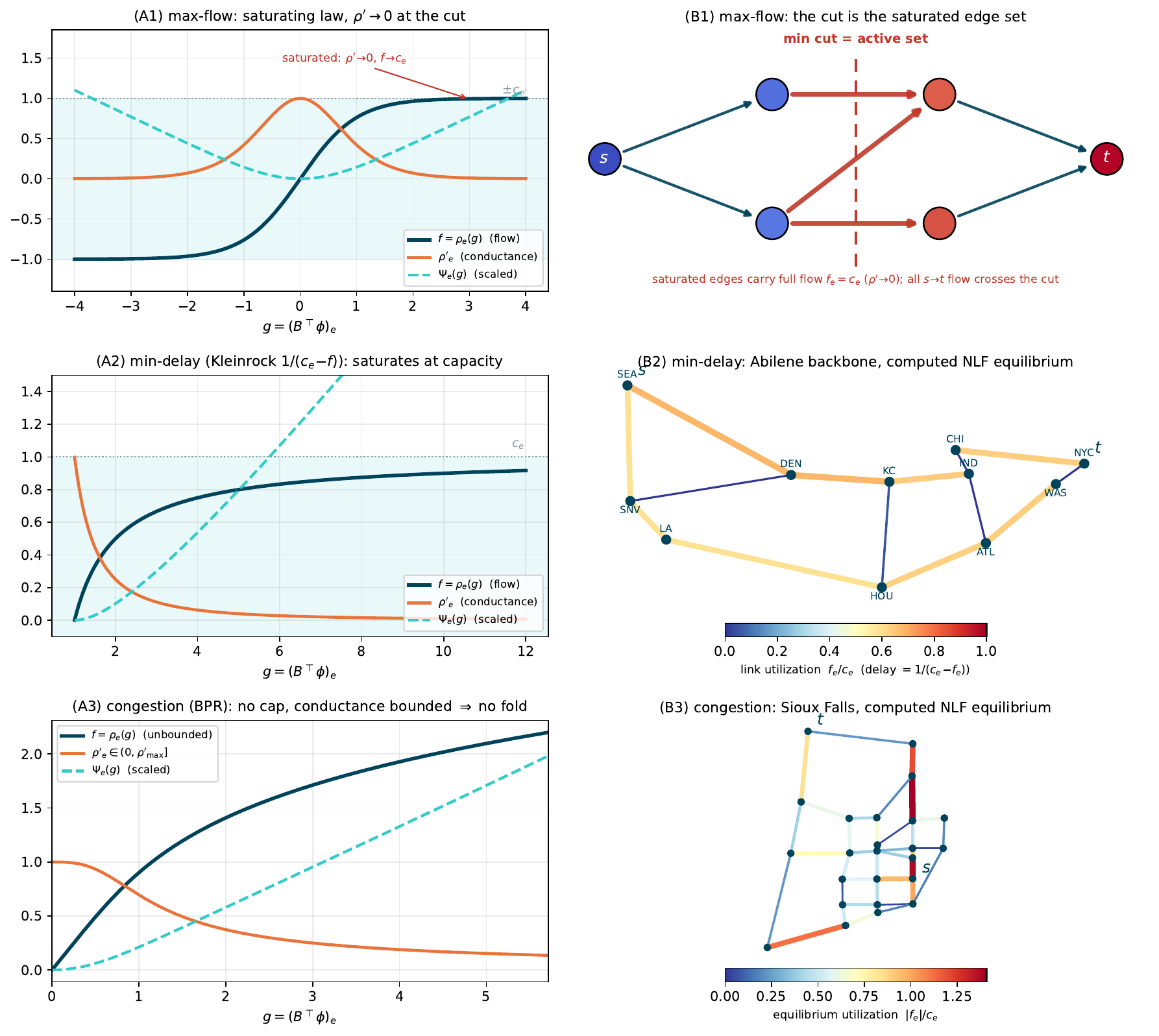}
\caption{The three instances of the framework (one row per line of Table~\ref{tab:instances});
column (A) the edge law $\rho_e$ (flow), its derivative $\rho'_e$ (the solution-dependent
conductance) and the energy $\Psi_e$; column (B) the behavior it induces.
\emph{Row 1, max-flow (tanh):} the law saturates at $\pm c_e$ and $\rho'_e\to0$; (B1) the saturated
edges are the min cut = the active set, and $\phi$ steps across them.
\emph{Row 2, minimum-delay routing (Kleinrock $1/(c_e\!-\!f)$):} the law saturates at capacity---the
fold class, handled by the continuation of \S\ref{sec:continuation}; (B2) the real application: the
Abilene/Internet2 backbone colored by the \emph{computed} NLF min-delay equilibrium for a
Seattle$\to$New York demand.
\emph{Row 3, congestion (BPR):} the law is unbounded and $\rho'_e$ is bounded away from zero---no
fold, pure cycling; (B3) the real application: the Sioux Falls road network colored by the
\emph{computed} NLF equilibrium utilization $|f_e|/c_e$ for a single source--sink demand
(\S\ref{sec:traffic}).}\label{fig:form}
\end{figure}

\subsection{Relation to prior Laplacian-paradigm formulations} Solving flow through graph-Laplacian
systems is the core of a productive line of work, but NLF's mechanism is fundamentally different. Electrical-flow methods~\cite{christiano-etal} and their $p$-Laplacian
descendants compute a sequence of \emph{linear} electrical flows $f=\mathrm{diag}(w)\,B^{\top}\phi$
---each a fixed weighted-Laplacian solve---and steer an \emph{outer} multiplicative-weights loop to a
$(1+\epsilon)$-approximate max-flow; IPM for (generalized) flow
\cite{daitch-spielman} likewise wrap a Laplacian solve in a barrier path-following loop. NLF instead places the nonlinearity in a \emph{constitutive law}: $f=\rho(B^{\top}\phi)$ is a
saturating function of the gradient, so the single energy of \eqref{eq:master} encodes the
capacity exactly and its feasibility limit is the exact min cut. The outer loop is then a short continuation through one scalar fold, not a reweighting or
barrier schedule. The Beckmann congestion instance \eqref{eq:primal} is already known
\cite{beckmann}; the novelty is the saturating-resistor reformulation of max-flow, and placing max-flow, congestion, and minimum-delay routing under one
nonlinear Laplacian solved via existing linear Laplacian algorithms.

\section{The NLF algorithm}\label{sec:algorithm}
The two problem classes need different outer loops. For congestion (no fold), $\alpha$ is a given input, so \eqref{eq:master} is solved \emph{directly} at that load by a handful of inexact chord-Newton steps (Algorithm~\ref{alg:fixed}), warm-started from $\phi=0$; no load ramp is needed---a single fixed-load solve converges across the corpus (\S\ref{sec:corpus}). For max-flow (fold), $\alpha$ is the unknown maximizer, and continuation from $\alpha=0$ in pseudo-arclength parametrization tracks the solution curve through the fold to $F^*$, each fixed-$(\phi,\alpha)$ subproblem again a handful of chord-Newton steps. The inner solver's setup---its factor or hierarchy---is built once, frozen across steps and re-formed only when an observed-convergence monitor signals it has drifted too far. The numbers of steps and inner solves are both bounded, implying empirical-$\Om$ complexity. We describe the linear engine
(\S\ref{sec:inner}), the fixed-load solve (\S\ref{sec:fas}), and the continuation (\S\ref{sec:continuation}).

\subsection{The inner Laplacian solve}\label{sec:inner}
The outer loop treats the inner solver as a black box that is required to do three things: (i)~it solves a
graph-Laplacian system $L\phi=b$ in near-linear $\Om$ work; (ii)~its accuracy is tunable, so it can match the $\eta\|r\|$ forcing tolerance  of the chord-Newton step (\S\ref{sec:fas}) rather
than to machine precision; and (iii)~its setup is separable from its solve, so the setup can be 
built once at a frozen linearization and reused across many continuation steps.

Two public near-$\Om$ solvers meet all three with different internals. \emph{Approximate
Cholesky}~\cite{approxchol} sets up a randomized incomplete sparse-Cholesky factor preconditioner and solves by preconditioned conjugate gradients to
the requested tolerance. \emph{LAMG+}~\cite{lamgplus} sets up an aggregation-multigrid
hierarchy (relaxation-based affinity aggregation, low-degree Schur elimination, a guarded coarse
operator) and solves by a multigrid cycle with min-residual recombination. We use approximate
Cholesky as the default as it is the faster default across the
corpus at equal robustness (\S\ref{sec:ablation}). LAMG+ is as robust, while its wall-clock constant
differs by graph class.

While both inner solvers yield the same overall nonlinear convergence (\S\ref{sec:corpus}), one convenient property of LAMG+'s relaxation-based aggregation is that it is \emph{cut-respecting}: as the cut conductances vanish they
decouple relaxation across the cut, so test vectors decorrelate and the affinity declines to aggregate
across it. Its measured Asymptotic Convergence Factor (ACF) of $\approx0.01$ thus holds through the continuation with
no fold-aware tuning (\S\ref{sec:maxflow}). This is a specific feature of LAMG+, not a framework requirement: adequate refresh of the frozen operator near the fold
suffices for any inner solver, and over-freezing degrades all of them equally. ($\alpha$, the
continuation, and the fold never enter the linear solve; the singular cut direction near $F^*$ is
deflated outside the engine, \S\ref{sec:continuation}.)

\subsection{The fixed-load solve: a frozen-setup inner solve}\label{sec:fas}
Eq.~\eqref{eq:master} with fixed $\alpha$ is solved by inexact chord-Newton iteration (Algorithm~\ref{alg:fixed}). At each step we evaluate the residual $r(\phi)=\alpha d-B\rho(B^{\top}\phi)$ and solve for the Newton direction $\delta\approx J_H^{-1}r$ using the inner Laplacian solver to tolerance $\eta\|r\|$, $\eta=0.05$~\cite{inexact-newton}. The Jacobian $J_H$ and the expensive setup hierarchy are only recalculated when chord Newton's linear convergence rate rises above $0.25$. This addresses the no-fold congestion class; we verified that $\phi=0$ is always a good initial guess there.
\begin{algorithm}[t]
\caption{NLF fixed-load solve}
\label{alg:fixed}
\begin{algorithmic}[1]
\Require load $\alpha$; initial guess $\phi$; inner-solver state $H$ (Cholesky factor or hierarchy, shared across calls); relative
residual tolerance $\mathrm{tol}$ (default $10^{-9}$); forcing parameter $\eta$ (default $0.05$)
\State $r \gets \alpha d - B\rho(B^{\top}\phi)$
\While{$\|r\| > \mathrm{tol}\cdot\alpha\,\|d\|$}
\State \textbf{refresh check:} if $H$ is empty, or the previous pass reduced $\|r\|$ by a factor
worse than $0.25$, or the line search stalled on a stale $H$, rebuild $H\gets$ inner setup of
$J(\phi)=B\,\mathrm{diag}(\rho'(B^{\top}\phi))B^{\top}$
\State $\delta \gets$ approximate solution of the frozen correction equation $J_H\,\delta = r$,
where $J_H=B\,\mathrm{diag}(\rho'(B^{\top}\phi_{\mathrm{build}}))B^{\top}$ is the linearization
$H$ was built from: inner solve(s) with $H$ starting from $\delta=0$, stopped when
$\|r - J_H\,\delta\| \le \eta\,\|r\|$
\State \textbf{line search:} for $\tau=1,\tfrac12,\tfrac14,\dots,\frac{1}{2^{49}}$ find the first $\tau$ with
$\big\|\alpha d - B\rho\big(B^{\top}(\phi+\tau\delta)\big)\big\| \le \|r\|$
\If{line search succeeded} $\phi\gets\phi+\tau\delta$
\ElsIf{$H$ was stale (built before this iteration)} rebuild $H$ from current $J(\phi)$;
recompute $\delta\gets$ inner solve of $J_H\,\delta=r$ from $\delta=0$; repeat line search
  \If{retry succeeded} $\phi\gets\phi+\tau\delta$
  \Else\ \textbf{break} \Comment{fresh-$H$ stall: no improving direction exists; exit}
  \EndIf
\Else\ \textbf{break} \Comment{fresh $H$ already in use; exit}
\EndIf
\State $r \gets \alpha d - B\rho(B^{\top}\phi)$
\EndWhile
\State \Return $\phi$
\end{algorithmic}
\end{algorithm}

\paragraph{The damped update} The line search (step 5) is standard backtracking on the residual
norm~\cite[Ch.~3]{nocedal-wright}. Because $\delta$ comes from a frozen linearization, a full step
can overshoot right after a load increase where the law's curvature is strongest; backtracking by
halving and accepting the first non-increasing residual guarantees monotone progress at negligible
cost. Near the solution the full step is always accepted (one trial, $\tau=1$).

\paragraph{When the setup is refreshed} Three conditions trigger: (i)~no setup yet; (ii)~per-step residual factor $>0.25$ (frozen $J$ has drifted); (iii)~line-search stall on a stale $H$---rebuild and retry once; a fresh-$H$ stall breaks the loop. The corpus median refresh count is small (often one). A lazy refresh is $1.5$--$4.7\times$ faster than rebuilding every step.

\paragraph{Deep hierarchies (LAMG+ engine only)} On graphs with ${\ge}20$-level hierarchies (country-scale roads, large FEM meshes) the dipole right-hand side excites the longest mode that the stock cycle under-treats. A \emph{grown} cycle index ($1.15\times$ per level, cap $\gamma_l\tau_l\le0.95$, each cycle stays $\Om$) is used, restarted if the factor exceeds $0.6$. This engages on $11$ of $2{,}003$ graphs ($0.5\%$) and restores textbook convergence.

\subsection{Continuation}\label{sec:continuation}
\begin{figure}[t]\centering
\begin{minipage}[b]{0.512\textwidth}\centering\includegraphics[width=\linewidth]{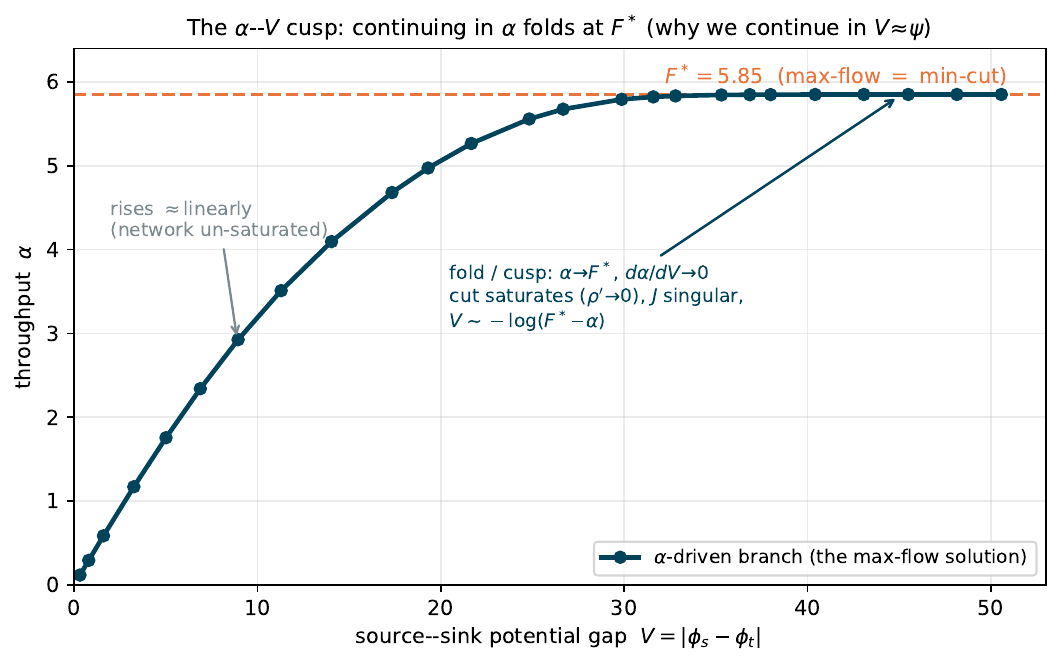}\\{\footnotesize(a)}\end{minipage}\hfill
\begin{minipage}[b]{0.447\textwidth}\centering\includegraphics[width=\linewidth]{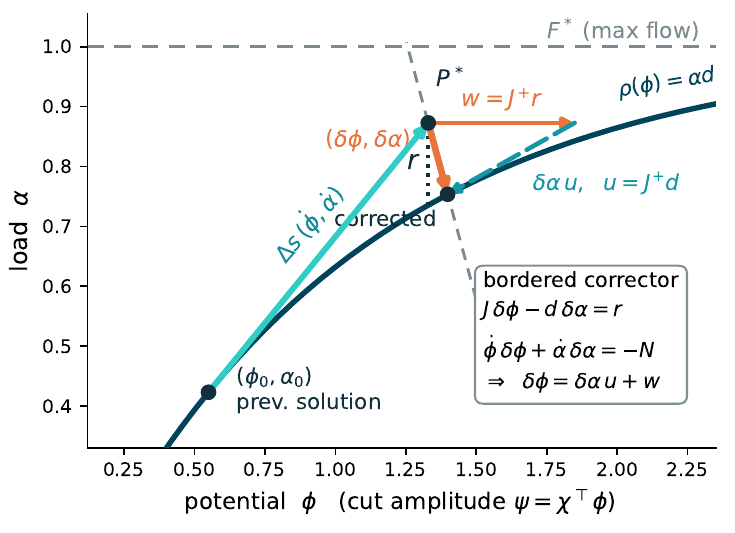}\\{\footnotesize(b)}\end{minipage}
\caption{\textbf{(a)}~The $\alpha$--$V$ fold (max-flow). As the load $\alpha$ is driven up, the
source--sink gap $V=|\phi_s-\phi_t|$ grows; $\alpha$ rises, then folds onto the min-cut value $F^*$ as
the cut saturates ($\rho'\to0$, $J$ singular) with $V\sim-\log(F^*-\alpha)$. Stepping $\alpha$ crowds
the fold; continuing in $V\approx\psi$ keeps steps uniform. Congestion laws have no such fold.
\textbf{(b)}~Pseudo-arclength predictor--corrector, drawn where $\phi,\alpha$ are both scalars so the
equilibria form a plane curve $\rho(\phi)=\alpha d$. From the previous solution the \emph{predictor}
advances $\Delta s$ along the unit tangent $(\dot\phi,\dot\alpha)$ to an off-curve point $P^*$; the
\emph{corrector} drives the equilibrium residual $r$ and the arclength residual $N$ to zero through the
bordered system \eqref{eq:bordered}. Its Newton step splits as $\delta\phi=\delta\alpha\,u+w$ (block
elimination \eqref{eq:deflate}) with $u=J^{+}d$ along the tangent and $w=J^{+}r$ cancelling $r$; near
the fold $u\propto\chi/\lambda_{\min}$ is large but is pinned by the scalar arclength row, so the inner
solves stay well-conditioned.}\label{fig:cusp}
\end{figure}

Only the saturating/fold class needs continuation. Here $\alpha$ is the unknown maximizer. A fixed-$\alpha$ solve is
feasible iff $\alpha<F^*$ (the min-cut value), and as $\alpha\uparrow F^*$ the solution curve
\emph{folds}---$\alpha(V)$ flattens onto the asymptote $F^*$ while the potentials grow without
bound (Fig.~\ref{fig:cusp}); at the fold $\kappa(J)\sim(F^*-\alpha)^{-1}$ and the source--sink gap
is $V=|\phi_s-\phi_t|\sim-\log(F^*-\alpha)$, so $dV/d\alpha\sim(F^*-\alpha)^{-1}\to\infty$.
Thus $\alpha$ is a bad continuation parameter: the feasible loads crowd at the fold and the step is forced to zero. The well-conditioned
variable is the \emph{cut-mode amplitude}
$\psi=\chi^{\top}\phi\sim-\log(F^*-\alpha)$,\footnote{Justification of both asymptotics. Near the fold
$\phi=V\chi+O(1)$, $\chi$ the source-side indicator, so every cut edge carries gradient
$g_e=V+O(1)$. The saturating tail of \eqref{eq:rho} is
$\rho_e(g)=c_e\tanh(g/c_e)=c_e-2c_e\,e^{-2g/c_e}+O(e^{-4g/c_e})$, and all $s$--$t$ flow crosses the
cut, so $\alpha=\sum_{e\in\mathrm{cut}}\rho_e(g_e)=F^*-2\sum_{e\in\mathrm{cut}}c_e\,
e^{-2V/c_e}\,(1+o(1))$. Hence $F^*-\alpha\asymp e^{-2V/\bar c}$, $\bar c$ the largest cut capacity,
i.e.\ $V=\tfrac{\bar c}{2}\log\tfrac{1}{F^*-\alpha}+O(1)$ and $\psi=\chi^{\top}\phi\propto V$. The
same expansion gives the cut conductances
$\rho'_e(g_e)=\mathrm{sech}^2(g_e/c_e)=4e^{-2V/c_e}(1+o(1))$, whence
$\lambda_{\min}(J)\le\sum_{\mathrm{cut}}\rho'_e\propto F^*-\alpha$ (here and below $\lambda_{\min}$
is the smallest \emph{nonzero} eigenvalue---$J$ being a graph Laplacian, its zero eigenvalue with
constant eigenvector $\mathbf 1$ is excluded, and $\kappa(J)$ is correspondingly the ratio over the
nonzero spectrum; the bound is the Rayleigh quotient at $\chi$,
which loads only the cut edges); for a \emph{simple} fold the next eigenvalue stays $\Theta(1)$ (the
deflated spectral gap) and $\lambda_{\max}(J)=\Theta(1)$, so the bound is two-sided and
$\kappa(J)\sim(F^*-\alpha)^{-1}$.}
$\chi$ the (a priori unknown) min-cut indicator. We realize this by pseudo-arclength
continuation~\cite{keller,allgower-georg} (Algorithm~\ref{alg:cont}), which never forms
$\chi$.

The equilibria of \eqref{eq:master} trace a smooth curve $\big(\phi(s),\alpha(s)\big)$ as the load is
varied; we parametrize it by its own arclength $s$ and write a dot for the derivative along the
curve, $\dot{(\cdot)}\equiv d(\cdot)/ds$, so the tangent $(\dot\phi,\dot\alpha)$ is a unit
vector, $\|\dot\phi\|^2+\dot\alpha^2=1$. Differentiating \eqref{eq:master} along $s$ gives:
\begin{equation}\label{eq:tangent}
J\dot\phi=\dot\alpha\,d\quad\Longrightarrow\quad \dot\phi=J^{+}d,\ \ \dot\alpha=1,\ \ \text{then rescale to }\|(\dot\phi,\dot\alpha)\|=1 .
\end{equation}
As $\alpha\to F^*$, $J^{+}d\propto\chi/\lambda_{\min}$ (the pseudoinverse $J^{+}$ acts on
$\mathbf 1^{\perp}$, so $\lambda_{\min}$ is again the smallest nonzero eigenvalue), so the tangent
rotates onto the cut mode by
itself---stepping in $s$ becomes stepping in $\psi$, the well-conditioned direction.

Each continuation step is a \emph{predictor} that advances a fixed arclength $\Delta s$ along the
tangent, $(\phi,\alpha)\mathrel{+}=\Delta s\,(\dot\phi,\dot\alpha)$, then a \emph{corrector} that
drives two residuals to zero: the \emph{equilibrium residual} $r=\alpha d-B\rho(B^{\top}\phi)$ of
\eqref{eq:master} (as in Algorithm~\ref{alg:fixed}) and the \emph{arclength residual}
\begin{equation}\label{eq:arclength}
N=\dot\phi^{\top}(\phi-\phi_0)+\dot\alpha\,(\alpha-\alpha_0)-\Delta s ,
\end{equation}
how far the point has advanced along the tangent past the previous solution $(\phi_0,\alpha_0)$, minus
the target $\Delta s$. Setting $r=0$ holds the point on the equilibrium curve; $N=0$ fixes how far
along it (hence
``pseudo''-arclength: the exact arclength sphere $\|\delta\phi\|^2+\delta\alpha^2=\Delta s^2$
is linearized to its tangent plane $N=0$.) One Newton step on $(r,N)=0$ is the bordered linear system
\begin{equation}\label{eq:bordered}
\begin{bmatrix} J & -d\\ \dot\phi^{\top} & \dot\alpha\end{bmatrix}
\begin{bmatrix}\delta\phi\\ \delta\alpha\end{bmatrix}=\begin{bmatrix}r\\ -N\end{bmatrix}.
\end{equation}
Block elimination of $\delta\alpha$ gives
\begin{equation}\label{eq:deflate}
\delta\alpha=-\,\frac{\dot\phi^{\top}w+N}{\dot\alpha+\dot\phi^{\top}u},\qquad
\delta\phi=\delta\alpha\,u+w ,
\end{equation}
where $u=J^{+}d$ and $w=J^{+}r$ are computed by two applications of the inner Laplacian solver.
Both $u,w\in\mathrm{range}(J^{+})=\chi^{\perp}$: the near-null cut component of the solution (the
$\chi$ direction) is pinned by the scalar arclength row \eqref{eq:arclength}, not by these solves. The
border keeps the system nonsingular through the fold ($\kappa=\mathcal{O}(10^3)$ while
$\kappa(J)\to\infty$), so the inner Laplacian solver stays fast (approximate Cholesky PCG and LAMG+ cycle ACF both $\approx0.01$).

\begin{algorithm}[t]
\caption{NLF arclength continuation (max-flow fold).}
\label{alg:cont}
\begin{algorithmic}[1]
\Require demand $d$; maximize $\alpha$ to find $F^*$ \Comment{congestion needs no continuation: Algorithm~\ref{alg:fixed} at the target load}
\State $\phi \gets 0$, $\ H \gets \emptyset$ \Comment{trivial solution at $\alpha=0$}
\State \emph{Seed a feasible load:}\ \ $\alpha_{\mathrm{hi}} \gets \sum_{e\sim s}c_e$ ($\ge F^*$, source-cut bound); set $\alpha\gets 0.1\,\alpha_{\mathrm{hi}}$ and halve until Algorithm~\ref{alg:fixed} converges (so $\alpha<F^*$)
\State \emph{Bracket $F^*$:}\ \ try load $g\alpha$ ($g{=}1.7$) with Algorithm~\ref{alg:fixed}; accept ($\alpha\gets g\alpha$) while it converges, stop at the first failure ($g\alpha\ge F^*$)
\State refresh $H$; tangent $(\dot\phi,\dot\alpha)$ by \eqref{eq:tangent}; initialize arc-step $\Delta s$
\While{$|\dot\alpha| > \mathrm{tol}$} \Comment{$\dot\alpha \to 0$ marks the fold, $\alpha\to F^*$}
  \State \emph{predictor:} $(\phi,\alpha)\mathrel{+}= \Delta s\,(\dot\phi,\dot\alpha)$;\ \ \emph{corrector:} iterate \eqref{eq:bordered} with line search, $H$ lazily refreshed (Algorithm~\ref{alg:fixed})
  \State on success grow $\Delta s$ and reorient the tangent \eqref{eq:tangent}; on a fresh-$H$ failure set $\Delta s\gets\Delta s/2$ and re-predict
\EndWhile
\State \Return $\phi$, $\ F^* = \alpha$
\end{algorithmic}
\end{algorithm}

The refresh policy of \S\ref{sec:fas} is woven into the corrector as a staleness safeguard (Algorithm~\ref{alg:cont}). Both $u$ and $w$ computations reuse one frozen hierarchy; since the cut conductances drain exponentially in the cut amplitude, that hierarchy ages faster as the fold sharpens, and the monitor rebuilds it more often — about every 4 arclength steps at $0.97 F^*$, every 2 at $0.995 F^*$. The predictor's tangent $\dot\phi = J^{+}d$, the single direction most sensitive to $J$'s drift, is instead recomputed on a fresh hierarchy once per accepted step. The fixed-load solution remains machine-precise as close to the fold as $\alpha = 0.995 F^*$ (a fixed-$\alpha$ equilibrium is well-posed only for $\alpha < F^*$), while the arclength continuation tracks through the fold to the exact F* without stalling (Table~\ref{tab:mf}).

\section{Transportation: undirected congestion equilibrium}\label{sec:traffic}
This is NLF's sharpest result: on convex congestion flow---which has no combinatorial competitor and
whose operator never folds---NLF is an empirically linear-time solver, verified against a
state-of-the-art interior-point method on real road-network topologies.

\paragraph{Scope: an undirected, single-commodity} Throughout this section the congestion
problem is the undirected, single-commodity Beckmann relaxation: flows are signed, the cost
depends on $|f_e|$, and Wardrop's path-nonnegativity is dropped (App.~\ref{app:instances}).
Our experiments establish two claims: (i) this undirected congestion
equilibrium is exactly \eqref{eq:master} for the inverse-BPR law, and (ii) a robust $\Om$
linear graph Laplacian inner solve makes it linear-time. The real road networks supply real topologies,
capacities, and BPR data exercised with a single source--sink demand; we cross-check NLF against an
independent high-accuracy solver on the identical convex program.

\subsection{The BPR congestion law}\label{sec:bpr}
On a congested road the travel time rises with the flow it carries. The field-standard model is the
Bureau of Public Roads (BPR) cost~\cite{bpr},
\begin{equation}\label{eq:bpr}
t_e(f)=t^0_e\Big(1+b\,(f/c_e)^{p}\Big),\qquad b=0.15,\ p=4 ,
\end{equation}
with $t^0_e$ the free-flow travel time and $c_e$ the practical capacity. The undirected Beckmann
congestion equilibrium (signed flows; Wardrop's path-nonnegativity relaxed, App.~\ref{app:instances})
is the convex program \eqref{eq:primal} with
$\Phi_e(f)=\int_0^f t_e=t^0_e f+\tfrac{b\,t^0_e}{p+1}f^{p+1}/c_e^{p}$; its edge law is the inverse
marginal cost $\rho_e=t_e^{-1}$, with conductance
$\rho'_e=1/t_e'(f)=\big(t^0_e b\,p\,(f/c_e)^{p-1}/c_e\big)^{-1}>0$ at every finite flow: the cost
rises but never caps the flow (the law does not saturate), so there is no feasibility limit and no fold.\footnote{We use a strictly convex, twice-differentiable form of
\eqref{eq:bpr} with the same free-flow and degree-$p$ congestion structure and real $(t^0_e,c_e,b,p)$
from the data; this regularizes the free-flow threshold so the dual problem is well posed. ``Practical
capacity'' sentinels coded as uncapacitated links are given a large finite capacity.}
The Jacobian $J=B\,\mathrm{diag}(\rho')B^{\top}$ is a well-conditioned Laplacian at every
bounded-flow iterate, and NLF runs the fixed-load Algorithm~\ref{alg:fixed}.

\subsection{No combinatorial competitor; an interior-point baseline}\label{sec:competitor}
No combinatorial algorithm addresses the convex congestion program \eqref{eq:primal}.
Frank--Wolfe~\cite{frankwolfe} converges only sublinearly on the Beckmann program: its
$\mathcal{O}(1/k)$ rate ($k$ = iteration count) drives the iteration count to even a $10^{-4}$ duality
gap up with problem size and cannot reach the high accuracy NLF attains---the classical slow-tail
limitation of the method on traffic assignment~\cite{boyles-tna}. The field's high-accuracy standard, the
bush-based family (Algorithm~B~\cite{bargera2002}, TAPAS~\cite{bargera2010}), solves the
\emph{directed} user equilibrium (nonnegative arc flows, one-sided cost)---a different convex program
from the undirected, signed-flow Beckmann relaxation NLF solves here (App.~\ref{app:instances}), so a
same-instance head-to-head is ill-posed. We therefore position bush methods qualitatively and
benchmark quantitatively against Ipopt~\cite{ipopt}, a mature IPM solver on the identical program. This is conservative: Ipopt attains the same $10^{-10}$, whereas
Frank--Wolfe and bush-based methods do not address the max-flow instance (\S\ref{sec:maxflow}).
Ipopt's per-step cost is a sparse-direct KKT factorization: cheap on planar graphs but
fill-in grows superlinearly as separators worsen---exactly the axis on which NLF is immune.

\subsection{Real road data: NLF and Ipopt}\label{sec:bprvalid}
We use the Transportation Networks for Research benchmark~\cite{tntp}, the field-standard
repository of real metropolitan road-network topologies with real capacities, free-flow times and BPR
parameters. Across instances spanning three decades of size, NLF and Ipopt reach the identical
equilibrium of the undirected single-commodity Beckmann program
(Table~\ref{tab:bprvalid}): the Beckmann objective agrees to the solver tolerance and the
equilibrium link flows to $\|f_{\mathrm{NLF}}-f_{\mathrm{Ipopt}}\|/\|f_{\mathrm{Ipopt}}\|\sim10^{-10}$,
in $\approx9$ Newton steps with no size trend.

\begin{table}[t]\centering\small
\caption{NLF (BPR congestion, near-linear inner) vs.\ Ipopt~\cite{tntp} on the \emph{undirected,
single-commodity} Beckmann congestion program over real road-network topologies. Both solvers reach
the \emph{same equilibrium of this program}; the objective matches to the tolerance, link flows to
$\sim10^{-10}$. This is the undirected net-flow relaxation, \emph{not} the field's directed
origin--destination user equilibrium (\S\ref{sec:traffic})---it is an exact cross-check of two solvers
on one identical convex program, so its objective value is not comparable to published directed-UE
results for these networks. NLF's step count shows no size trend. (Two further instances,
Philadelphia and Berlin-Center, also match to $10^{-10}$ but need more steps under extreme capacity
heterogeneity; omitted.)}\label{tab:bprvalid}
\begin{tabular}{lrrrrr}
\toprule
network & $n$ & $m$ & NLF obj.\ & $\|\Delta f\|/\|f\|$ & steps\\
\midrule
Sioux Falls       & 24     & 38     & $784.96$    & $3{\cdot}10^{-9}$  & 7\\
Anaheim           & 416    & 634    & $1600.92$   & $6{\cdot}10^{-12}$ & 9\\
Chicago-Sketch    & 933    & 1\,475 & $8178.62$   & $3{\cdot}10^{-11}$ & 9\\
Austin            & 7\,388 & 10\,591& $16615.74$  & $4{\cdot}10^{-12}$ & 10\\
Chicago-regional  & 12\,979& 20\,627& $281.15$    & $2{\cdot}10^{-11}$ & 9\\
Sydney            & 32\,956& 38\,787& $1198.98$   & $2{\cdot}10^{-11}$ & 9\\
\bottomrule
\end{tabular}
\end{table}

\subsection{Linear scaling, and the direct-vs-iterative crossover}\label{sec:bprscale}
\begin{figure}[t]\centering
\includegraphics[width=\textwidth]{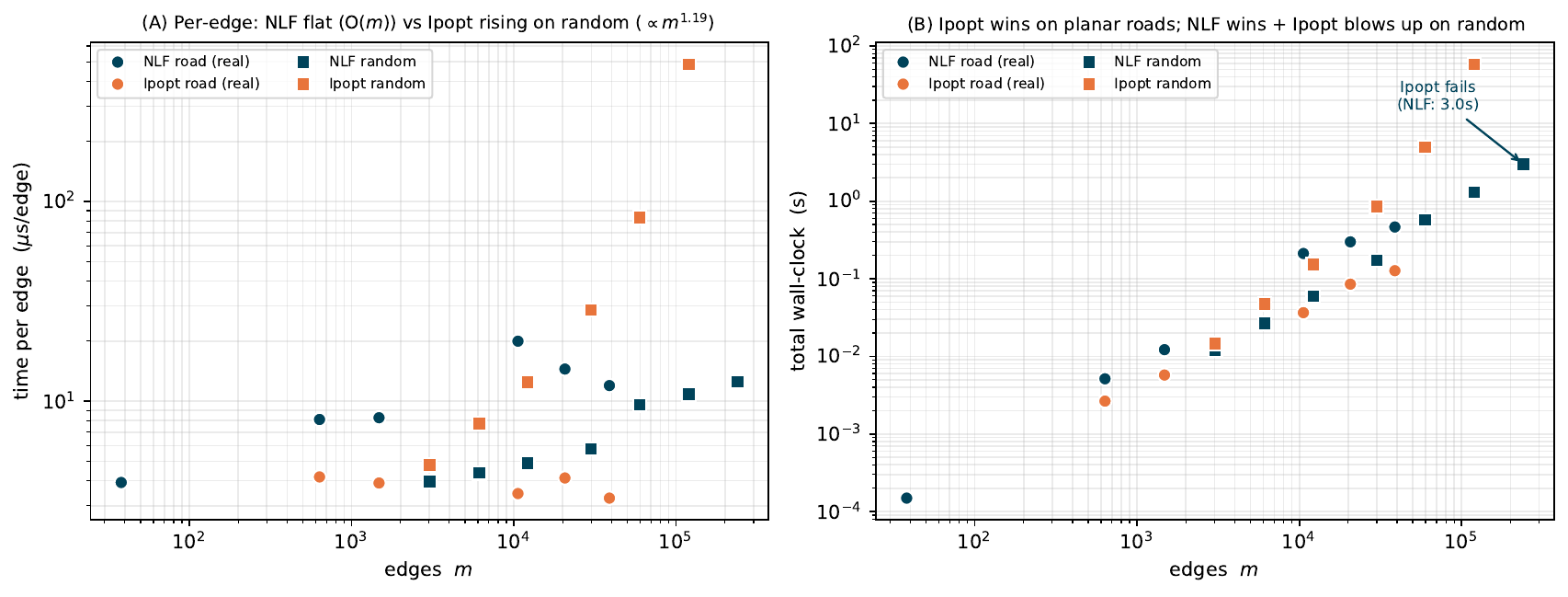}
\caption{Congestion (BPR) equilibrium, NLF vs.\ Ipopt, on real road networks (circles) and synthetic
poorly-separable graphs (squares). (A)~\emph{Per-edge} time: NLF is flat (near-$\Om$) on both
families; Ipopt is flat on roads but rises steeply on random graphs ($\propto m^{1.2}$ per edge).
(B)~Total wall-clock: Ipopt wins on planar roads (cheap separators), is overtaken on random graphs by
a margin that widens to $45\times$, and fails (time limit) at $m=2.4\times10^5$ where NLF finishes
in $3.0$\,s.}\label{fig:traffic}
\end{figure}

Figure~\ref{fig:traffic} and Table~\ref{tab:bprscale} time both solvers on the real road networks
(planar, good separators) and on synthetic Erd\H{o}s--R\'enyi graphs carrying BPR parameters (poorly
separable), to $2.4\times10^5$ edges. NLF's per-edge cost is flat on both families---empirically
near-linear (the full-corpus fit of \S\ref{sec:corpus} gives $t\propto m^{1.02}$); neither the
step count ($\approx9$) nor the inner linear solver ACF sees the graph's separability, while Ipopt does. On the real road networks its sparse-direct core is
near-linear and beats NLF by $\approx3\times$; planar graphs are where direct methods
excel. On poorly-separable graphs the KKT fill-in explodes---$t\propto m^{2.2}$---so the gap inverts
and widens with size: NLF is $5\times$ faster at $m=3.0\times10^4$, $45\times$ at $m=1.2\times10^5$,
and at $m=2.4\times10^5$ Ipopt does not finish within its $120$\,s CPU budget while NLF finishes in
$3$\,s.
\begin{table}[t]\centering\small
\caption{Scaling on poorly-separable (Erd\H{o}s--R\'enyi) graphs with BPR cost: NLF
($\Om$ near-linear inner) vs.\ Ipopt 3.14.19 (interior-point, MUMPS~5.9.0 sparse-direct KKT). Both reach the
same optimum of the undirected congestion program. NLF's step count shows no size trend; Ipopt is
superlinear and exceeds the $120$\,s CPU budget on the largest instance (``---'').
Speed-up $=t_{\mathrm{Ipopt}}/t_{\mathrm{NLF}}$.}\label{tab:bprscale}
\begin{tabular}{rrrrr}
\toprule
$m$ & NLF (s) & steps & Ipopt (s) & speed-up\\
\midrule
12\,175  & 0.060 & 9 & 0.152 & $2.5\times$\\
29\,939  & 0.173 & 9 & 0.856 & $4.9\times$\\
59\,732  & 0.573 & 9 & 4.97  & $8.7\times$\\
119\,460 & 1.30  & 9 & 58.2  & $45\times$\\
239\,842 & 3.01  & 9 & \multicolumn{1}{c}{--- ($>120$\,s)} & \multicolumn{1}{c}{---}\\
\bottomrule
\end{tabular}
\end{table}

\subsection{Robustness across the graph corpus}\label{sec:corpus}
The road networks above exercise NLF on real flow topologies; we now show it is robust across graph classes in full. We impose a single-commodity BPR congestion instance (random
capacities and free-flow times, a far source--sink demand) on every graph of a $2{,}003$ real-world
SuiteSparse collection, each reduced to its largest connected
component. This includes spanning structured grids, finite-element and circuit meshes, web,
social, road, and citation networks with up to $1.8\times10^7$ edges.

\begin{figure}[t]\centering
\includegraphics[width=\textwidth]{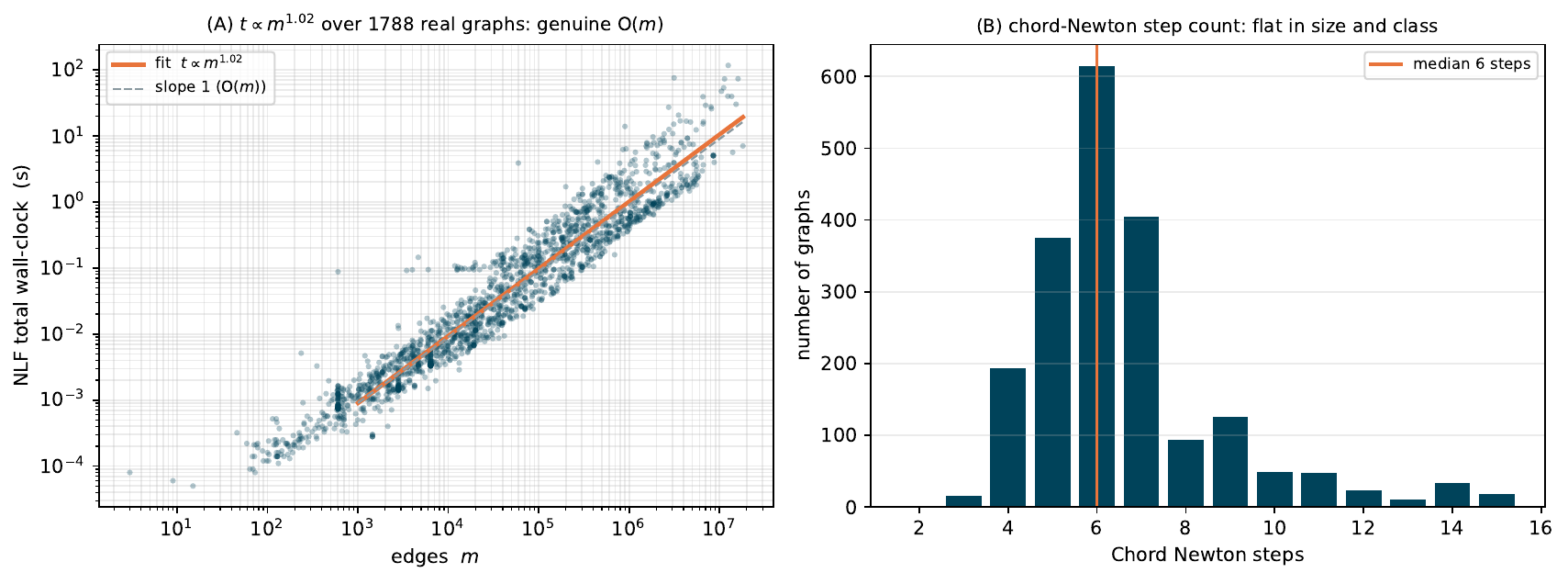}
\caption{NLF over the full real-world SuiteSparse corpus ($2{,}003$ graphs, all classes,
$10^2$--$1.8\times10^7$ edges, $100\%$ converged). (A)~total wall-clock vs.\ $m$ (log--log): the fit
is $t\propto m^{1.02}$ over five decades (slope-$1$ reference dashed), empirically near-linear,
per-edge cost flat at ${\approx}0.86~\mu$s/edge median. (B)~Newton step count: a flat distribution
(median $6$, range $2$--$15$).}\label{fig:corpus}
\end{figure}

Under this protocol (one BPR instance per graph, random capacities, a single far source--sink
demand, stock settings) NLF converged on all $2{,}003$ graphs to $10^{-9}$ residual tolerance, with a Newton step count median of $6$ and no size trend (Fig.~\ref{fig:corpus}B). The wall-clock fits \textbf{$t\propto m^{1.02}$} over the $1{,}788$ graphs above $10^3$ edges (Fig.~\ref{fig:corpus}A),
per-edge cost flat at ${\approx}0.86~\mu$s/edge median ($95$th percentile $3.7~\mu$s/edge). Three giants
confirm the scaling beyond the sweep's cap, to \texttt{hollywood-2009} at $5.6\times10^7$ edges.

\subsection{Comparison with other convex optimizers.}\label{sec:frontier}
We compare three solver paradigms for the same convex BPR program: near-linear (NLF), sparse-direct
(Ipopt), and matrix-free first-order (L-BFGS on the dual \eqref{eq:dual}; on well-conditioned graphs
L-BFGS also scales near-linearly and beats Ipopt, so that win is a fill-in effect any matrix-free
method inherits). We ran all three on the $1{,}669$ corpus graphs with $n\le300{,}000$ (the subset where Ipopt's sparse-direct KKT does not run out of memory on poorly-separable instances), each competitor capped at $5\times$ NLF's wall-clock (Table~\ref{tab:robustness}).
NLF converged on all $1{,}669$; first-order failed on $481$ (ill-conditioning,
$\kappa\sim N$) and interior-point on $143$ (fill-in). On the $90$ graphs hard on both axes,
neither competitor finished within budget. Where they converge, NLF is median
$2.6\times$/$4.2\times$ faster than Ipopt/L-BFGS and fastest-or-only on $85\%$.
Timing protocol: all three solvers share one global precompilation pass, and NLF is then timed on a
warm per-graph call; the competitors on their first per-graph call. This per-graph warm-up favours NLF
by a small-graph-weighted margin that is negligible at scale, where NLF's advantage is the genuine
avoidance of direct fill-in (\S\ref{sec:bprscale}); the ``within budget'' counts use the $5\times$-NLF
wall-clock cap and are a wall-clock criterion, not a claim the others cannot converge given more time.

\begin{table}[t]\centering\small
\caption{Three-paradigm robustness frontier across the SuiteSparse corpus ($1{,}669$ graphs, undirected
BPR congestion---the same convex program for all). NLF (multigrid-Newton) vs.\ Ipopt
(interior-point/sparse-direct KKT) vs.\ L-BFGS-on-the-dual (matrix-free first-order); each competitor
capped at $5\times$ the measured NLF wall-clock ($1$\,s floor, $120$\,s ceiling). NLF converges on every
graph and is the only paradigm that stays within budget on both axes. On $90$ graphs neither competitor
finishes within the $5\times$NLF budget, and Ipopt is additionally too large to run on $2$ more---a
union of $92$ on which NLF alone finishes (``within budget'' is a wall-clock criterion, not a claim the
others cannot converge given more time).
}\label{tab:robustness}
\begin{tabular}{p{0.34\columnwidth} r p{0.44\columnwidth}}
\toprule
regime & graphs & outcome\\
\midrule
both competitors converge (well-conditioned) & $1{,}137$ & NLF $\le$ both (median $2.6\times$/$4.2\times$ faster than Ipopt/L-BFGS)\\
ill-conditioned (first-order degrades) & $389$ & L-BFGS exceeds $5\times$ NLF; NLF, Ipopt converge\\
poorly-separable (direct fills in) & $51$ & Ipopt over budget (fill-in); NLF, L-BFGS converge\\
\textbf{both over budget $\Rightarrow$ union} & \textbf{$92$} & \textbf{NLF is the only solver that finishes}\\
\midrule
\textbf{total} & \textbf{$1{,}669$} & \textbf{NLF $100\%$ converged; fastest-or-only on $85\%$}\\
\bottomrule
\end{tabular}
\end{table}

\subsection{Inner-solver ablation: which near-linear engine}\label{sec:ablation}
The framework can use any near-linear inner solver (\S\ref{sec:inner}); we settle the
default empirically. Holding NLF's outer method fixed (load-continuation chord-Newton, lazy refresh),
we swap only the inner Laplacian solve and rerun the congestion corpus of
\S\ref{sec:corpus} (one BPR instance each; over $1{,}800$ graphs except the very largest graphs). The candidates are the two
public $\Om$ solvers of \S\ref{sec:inner}, approximate Cholesky and LAMG+; we exclude sparse-direct,
which is not near-linear---its separability crossover against an $\Om$ inner is the scaling result of
\S\ref{sec:bprscale}, not an inner we would ship.

\emph{Both converge on every graph, to the identical equilibrium}: the dual energy agrees to
floating-point precision and the link flows to $\le2\times10^{-8}$ relative error. So robustness is therefore not a discriminator. Speed splits by graph class, both ways: approximate Cholesky is faster on $81\%$ of
graphs (median $1.8\times$, up to $43\times$---its cheap factor wins on low-treewidth mesh and
finite-element graphs) and LAMG+ on the other $19\%$ (up to $24\times$---its hierarchy wins on the
poorly-separable and structural graphs); the overall median is $1.6\times$ in approximate Cholesky's favour, with
$39\%$ of graphs within $1.5\times$ either way (Table~\ref{tab:innerswap}). We therefore adopt the faster
approximate Cholesky as the default and keep LAMG+ as a
robust alternative: it wins wall-clock on the poorly-separable and structural minority, and at the
saturating fold its cut-respecting aggregation (\S\ref{sec:inner}) gives a structurally clean
near-singular solve with no fold-aware tuning, though approximate Cholesky's constant is often smaller
even there (\S\ref{sec:innerfold}). 

Residual tolerance criterion sensitivity is mild (Table~\ref{tab:paramablation}): freezing the factor more
aggressively trades a few extra outer steps for fewer setups until it over-freezes, and a looser inner
forcing $\eta$ costs outer steps; the stock $(\text{refresh}=0.25,\ \eta=0.05)$ sits in the flat
interior.

\begin{table}[t]\centering\small
\caption{Inner-solver ablation inside NLF's BPR congestion solve, outer method fixed: LAMG+ vs.\
approximate Cholesky (both $\Om$). Wall-clock seconds and outer (chord-Newton) steps, to a $10^{-7}$
residual; \emph{both return the identical equilibrium} (energy $E$ to the digits shown; flows to
$\Delta f$ relative). $t_{\mathrm{L}}/t_{\mathrm{A}}{>}1$ means approximate Cholesky is faster. A
representative cross-class subset spanning both speed directions; the corpus-wide summary is in the
text and the full per-graph output in the released artifact.}\label{tab:innerswap}
\resizebox{\columnwidth}{!}{%
\begin{tabular}{llrrrrr}
\toprule
graph & class & $m$ & energy $E$ & LAMG+ (stp) & approxChol (stp) & $t_{\mathrm{L}}/t_{\mathrm{A}}$\\
\midrule
\texttt{3elt}          & 2D mesh     & 13{,}722  & $-0.8939$ & $0.050$ (16) & $0.015$ (21) & $3.3$\\
\texttt{airfoil1}      & 2D FEM      & 12{,}289  & $-0.9866$ & $0.048$ (18) & $0.015$ (21) & $3.3$\\
\texttt{bodyy5}        & aniso.\ FEM & 55{,}346  & $-0.9045$ & $0.56$ (18)  & $0.070$ (21) & $8.0$\\
\texttt{bcsstk38}      & structural  & 173{,}714 & $-0.0432$ & $0.14$ (15)  & $0.22$ (21)  & $0.6$\\
\texttt{G2\_circuit}   & circuit     & 288{,}286 & $-1.9134$ & $2.63$ (21)  & $0.58$ (21)  & $4.5$\\
\texttt{fe\_ocean}     & FEM mesh    & 409{,}593 & $-1.0609$ & $3.84$ (18)  & $0.96$ (21)  & $4.0$\\
\texttt{p2p-Gnutella04}& P2P         & 39{,}994  & $-3.1625$ & $0.14$ (24)  & $0.16$ (25)  & $0.9$\\
\texttt{as-caida}      & AS/comm     & 53{,}381  & $-6.6601$ & $0.083$ (20) & $0.17$ (23)  & $0.5$\\
\texttt{web-NotreDame} & web         & 1{,}090{,}108 & $-25.957$ & $4.51$ (27) & $3.69$ (29) & $1.2$\\
\bottomrule
\end{tabular}}
\end{table}

\begin{table}[t]\centering\small
\caption{Outer-loop parameter sensitivity (approximate-Cholesky inner): outer (chord-Newton) steps to
a $10^{-9}$ residual on representative congestion instances. \emph{Refresh}---rebuild the frozen factor
when a step's residual reduction is worse than this threshold ($\eta{=}0.05$ fixed): flat in the
interior, a few extra steps only when over-frozen (\texttt{as-caida}). \emph{Forcing} $\eta$---the
inner relative tolerance (refresh${=}0.25$ fixed): a looser inner trades directly into more outer
steps. The stock $(\text{refresh},\eta){=}(0.25,0.05)$ (bold) sits in the flat interior and is
near time-optimal.}\label{tab:paramablation}
\begin{tabular}{lrrrr@{\qquad}rrr}
\toprule
 & \multicolumn{4}{c}{refresh ($\eta{=}0.05$)} & \multicolumn{3}{c}{forcing $\eta$ (refresh${=}0.25$)}\\
\cmidrule(lr){2-5}\cmidrule(l){6-8}
graph & $0.10$ & $\mathbf{0.25}$ & $0.50$ & $0.80$ & $0.01$ & $\mathbf{0.05}$ & $0.20$\\
\midrule
\texttt{airfoil1}      & 18 & \textbf{18} & 18 & 18 & 14 & \textbf{18} & 33\\
\texttt{as-caida}      & 21 & \textbf{20} & 29 & 29 & 17 & \textbf{19} & 25\\
\texttt{bodyy5}        & 18 & \textbf{21} & 18 & 19 & 15 & \textbf{18} & 36\\
\texttt{delaunay\_n16} & 18 & \textbf{18} & 21 & 18 & 15 & \textbf{18} & 33\\
\bottomrule
\end{tabular}
\end{table}

\subsection{Cost as a function of the target accuracy}\label{sec:accuracy}
The dependence on the requested solution accuracy $\varepsilon$ is logarithmic. Each inner solve to
relative residual $\delta$ costs $\mathcal{O}(m\log(1/\delta))$ (a bounded inner contraction gives
$\mathcal{O}(\log(1/\delta))$ inner iterations). The outer iteration is a frozen chord-Newton with
inexact-Newton forcing $\eta_k$ tracking the nonlinear residual $\|r_k\|$ (\S\ref{sec:inner}); its
per-step cycle counts form a series dominated by the final solve (where
$\|r_k\|\!\to\!\varepsilon$), the earlier looser steps adding only a bounded geometric tail, while the
continuation and globalization run at a fixed, $\varepsilon$-independent path tolerance. Each
contribution is thus a bounded multiple of a single linear solve, so the total cost to accuracy
$\varepsilon$ is, empirically,
\begin{equation}\label{eq:eps}
\mathcal{O}\!\big(m\,\log(1/\varepsilon)\big).
\end{equation}
The constant is the whole-nonlinear-solve-to-one-linear-solve ratio reported below ($2$--$4$), and it
stays bounded because the frozen chord-Newton contraction stays bounded away from $1$: on the no-fold
congestion class directly (\S\ref{sec:traffic}), and through the fold because the continuation
deflates the diverging cut direction and keeps the bordered corrector uniformly conditioned
(\S\ref{sec:continuation}); empirically the contraction does not degrade at $F^*$. We verify this by
sweeping $\varepsilon$ from $10^{-2}$ to $10^{-12}$ on one graph per class (web, social, road,
finite-element): the total inner-solve count grows \emph{linearly} in $\log(1/\varepsilon)$, at
$0.24$--$0.60$ per decade of accuracy, confirming \eqref{eq:eps}.

\paragraph{The nonlinear/linear cost ratio is $\mathcal{O}(1)$} Brandt's yardstick says: solving the nonlinear problem should cost only
a small, size-independent multiple of solving one linear problem of the same
kind~\cite[\S8.3]{brandt-guide}. We measure exactly this: the full congestion solve (load
continuation, all steps, all cycles, lazy refresh) against a single inner setup-and-solve of
the linearization at $\alpha=0$ (the linear resistor network $B\,\mathrm{diag}(\rho'(0))B^{\top}\phi=d$)
on the same graph. Across the classes the ratio is
\textbf{$1.9$--$4.4$}---grid $3.1$, finite-element $2.5$, road $1.9$, social $3.8$, web ($1.9\times10^6$
edges) $4.4$---flat in size and class.

\subsection{How accurately to solve each corrector step at the fold}\label{sec:contin-ablation}
At the fold, continuation is mandatory ($J$ singular at $F^*$), but each arclength corrector need only
stay on the solution path~\cite[\S8.3.2]{brandt-guide}---a loose per-step tolerance reaches the same
min cut in $2$--$7\times$ fewer corrector solves and hierarchy rebuilds, while a single corrector step
per arclength step overshoots $F^*$.

\section{Maximum flow}\label{sec:maxflow}
The saturating instance \eqref{eq:rho} recovers the exact combinatorial maximum flow. Read
$f=\rho(B^{\top}\phi)$ as an edge-flow vector: $\rho$'s range is the capacity box $|f_e|\le c_e$, and
$Bf=\alpha d$ is conservation, so any $\phi$ solving \eqref{eq:master} is a capacity-feasible,
conserved flow of value $\alpha$. The recession slope of $\Psi_e$ is $c_e$, so $E$ has a finite
minimizer iff $\alpha<F^*$; as $\alpha\uparrow F^*$ the min-cut edges saturate ($\rho'_e\to0$ there
only) and $\phi$ develops a cut-indicator step. By max-flow/min-cut duality the largest feasible
$\alpha$ is $F^*$, the min-cut capacity, and $\phi$ is the LP dual. Because $f=\rho(B^{\top}\phi)$ is a
\emph{nonlinear} function of the gradient, it realizes the true max-flow; a \emph{linear} resistor
network computes only the sub-maximal electrical flow ($0.02$--$0.20\,F^*$ on heterogeneous
capacities, except single-bottleneck instances where the two coincide). The nonlinearity is essential.

\begin{figure}[t]\centering
\includegraphics[width=\textwidth]{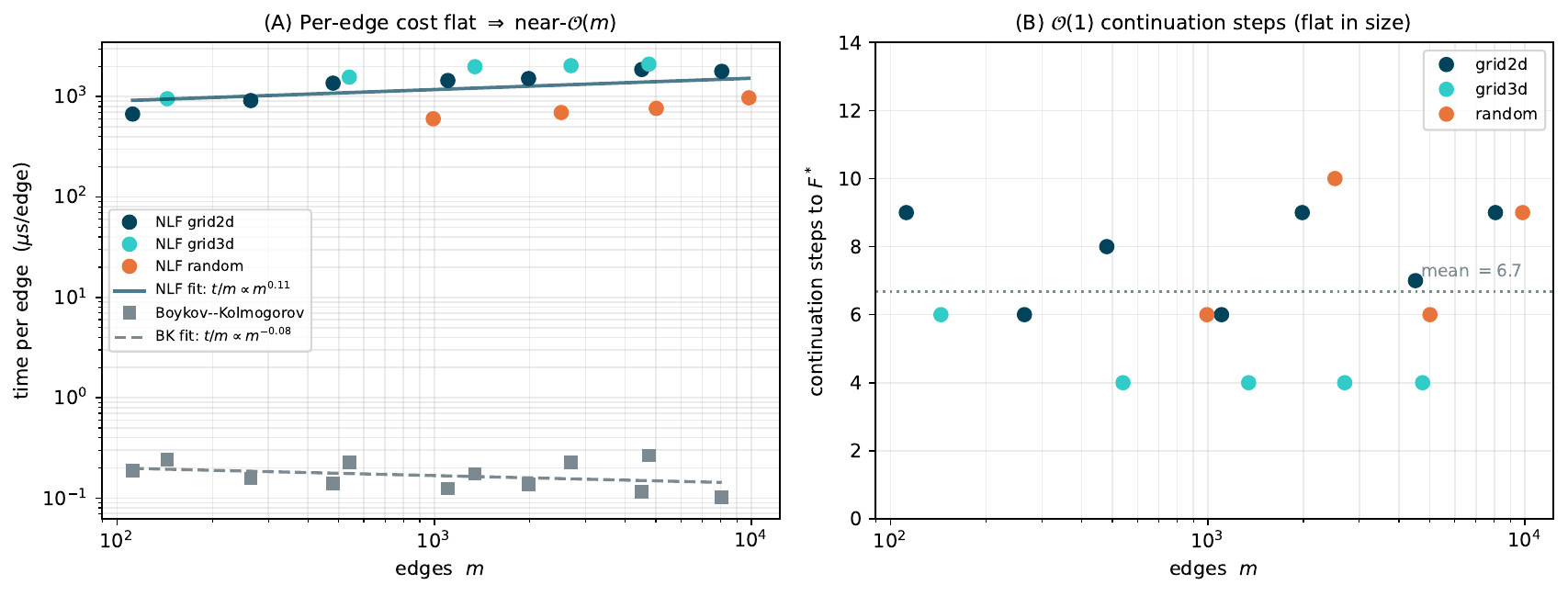}
\caption{Max-flow scaling (whole algorithm, total wall-clock), on 2D/3D grids and random graphs to
${\sim}10^4$ edges, all exact at $F^*$. (A)~\emph{Per-edge} cost: nearly flat
($t\propto m^{1.00}$, empirically $\Om$) at $\approx1.4$\,ms/edge, the constant set by the per-step
hierarchy rebuild; Boykov--Kolmogorov sits far below at $\approx0.17~\mu$s/edge (calibration only).
(B)~The continuation is $\mathcal{O}(1)$: pseudo-arclength steps are flat in size (mean
$6.7$).}\label{fig:scaling}
\end{figure}

\paragraph{Cut-respecting coarsening} As $\alpha\uparrow F^*$ the linearization $J$ becomes singular
along a single near-null direction---the min-cut indicator---so $\kappa(J)\to\infty$. We assume a
\emph{simple fold}: the min cut is unique, so the limiting near-null space is one-dimensional (a
degenerate cut of multiplicity $k$, e.g.\ in symmetric unweighted graphs, would make it $k$-dimensional
and require a rank-$k$ block border in place of the scalar one below). The continuation
of \S\ref{sec:continuation} \emph{deflates} that one direction to a scalar outside the engine (the
bordered/pseudo-arclength system); what the inner solver sees is the deflated complement
(Eq.~\ref{eq:deflate}). \emph{The
deflation is the enabling property}---it is what keeps the inner operator well-conditioned through the
fold, and it is a property of the framework, not of any one inner solver (swapping multigrid for a
sparse-direct or another near-linear Laplacian solve reaches the identical $F^*$;
\S\ref{sec:innerfold}).
Multigrid's relaxation-based affinity is a natural fit on the complement---since the forming cut is the
zero-conductance set of $J$ (\S\ref{sec:formulation}), the affinity declines to aggregate across it, so
the coarse hierarchy mirrors the cut without being told where it is---but it is helpful, not required.
The measured multigrid factor \emph{on the deflated operator} is $\mu\approx0.01$, uniformly in
$\alpha/F^*$ down to $0.999\,F^*$ (Fig.~\ref{fig:scaling}); the uniformity is a property of the
deflation, not of aggregation alone---on the undeflated $J$, $\mu$ degrades with $\kappa(J)$ as
expected.

\paragraph{Exactness and scaling} NLF returns the exact max-flow (Table~\ref{tab:mf}), checked
against $F^*$ from a linear program over free flows and from push--relabel (which agree to $10^{-4}$),
on every graph, to solver tolerance. The pseudo-arclength step count is flat in size
(Fig.~\ref{fig:scaling}B, mean $6.7$) and the total wall-clock is empirically $\Om$: $t\propto m^{1.00}$
over the max-flow sweep (to ${\sim}10^4$ edges, Fig.~\ref{fig:scaling}A). \emph{Calibration, not a contest:}
Boykov--Kolmogorov, on its target vision grids, runs at $\approx0.17~\mu$s/edge
(Fig.~\ref{fig:scaling}A), far below NLF's constant---\textbf{NLF is not meant to be a faster
combinatorial max-flow solver}. Its value is being $\Om$, returning the cut potential (LP dual), the
full flow-vs-load curve, and a smooth interior flow---the objects an outer method manipulates when
max-flow is an inner step---one instance of a framework that, on the congestion class, has \emph{no}
combinatorial competitor at all.

\begin{table}[t]\centering\small
\caption{NLF recovers the true max-flow exactly; the gradient/electrical relaxation falls far short.
$F^*$: LP over free flows and push--relabel agree to $10^{-4}$; \texttt{washington} and \texttt{genrmf}
are DIMACS-challenge generators~\cite{dimacs}.}\label{tab:mf}
\begin{tabular}{lrrcc}
\toprule
graph & $n$ & $F^*$ & grad$/F^*$ & NLF$/F^*$\\
\midrule
grid2d/$16$ & 256 & 4.973   & 0.20 & 1.0000\\
grid3d/$6$  & 216 & 17.478  & 0.11 & 1.0000\\
washington  & 66  & 500.745 & 0.02 & 1.0000\\
genrmf      & 64  & 16.000  & 0.00 & 1.0000\\
bottleneck  & 40  & 1.000   & 1.00 & 1.0000\\
\bottomrule
\end{tabular}
\end{table}

\subsection{Inner solver ablation}\label{sec:innerfold}
The deflation, not the choice of inner solver, is what carries the continuation through the fold
(\S\ref{sec:maxflow}); the inner Laplacian solve is therefore a replaceable module. We check this
directly: holding NLF's outer method fixed (chord-Newton with the deflated pseudo-arclength
continuation to $F^*$), we swap only the inner solver---multigrid (LAMG+), approximate
Cholesky~\cite{approxchol}, an exact sparse Cholesky factorization, and a diagonally preconditioned
CG---on max-flow instances spanning FEM, structural, and circuit graphs (Table~\ref{tab:innerfold}).

The result is modularity, not a winner: every inner solver that converges returns the identical
$F^*$ to solver tolerance, so the continuation is genuinely solver-agnostic. The distinction that
matters is near-linear versus not: the exact factorization is fastest where it fits but exhausts
memory past ${\sim}150$k poorly-separable nodes (\texttt{G2\_circuit}, \texttt{scircuit}), exactly the
regime where a near-linear inner solve is mandatory, while the diagonal-PCG baseline---no coarse
correction---is typically up to two orders of magnitude slower and the least scalable. Among the
near-linear solvers, approximate Cholesky and LAMG+ are interchangeable in robustness; approximate
Cholesky often carries the smaller constant. The single
stiff case (\texttt{bcsstk38}) is instructive: there \emph{all} iterative inner solves stall and only
the exact factorization reaches $F^*$---the deflation removes the \emph{global} cut singularity but not the residual ill-conditioning of the deflated system, for which an exact inner solve is still required.

\begin{table}[t]\centering\small
\caption{Inner-solver swap inside NLF's max-flow continuation to $F^*$ (wall-clock seconds; ``stall''
$=$ did not converge within $120$ steps; ``OOM'' $=$ exact factorization size-guarded at $150$k
nodes). Every converging solver returns the \emph{identical} $F^*$, so the continuation is
solver-agnostic; the operative choice is near-linear versus not. Representative subset; the bake-off
script and its outputs are in the released artifact. On \texttt{bcsstk38} the stalled iterative solvers were arrested at
$37.0$, far from the true $F^*\!=\!45.6$ the direct solve reaches.}\label{tab:innerfold}
\resizebox{\columnwidth}{!}{%
\begin{tabular}{llrrrrrr}
\toprule
graph & class & $n$ & $F^*$ & LAMG+ & approxChol & diagPCG & direct\\
\midrule
\texttt{grid2}       & 2D FEM      & 3{,}296   & 3.003  & 11.4  & 0.75  & 12.3  & 0.24\\
\texttt{airfoil1}    & 2D FEM      & 4{,}253   & 6.695  & 29.3  & 1.30  & 20.6  & 0.40\\
\texttt{bodyy5}      & aniso.\ FEM & 18{,}589  & 3.325  & 307   & 9.5   & 214   & 2.0\\
\texttt{bcsstk38}    & structural  & 8{,}032   & 45.61  & stall & stall & stall & 4.0\\
\texttt{circuit\_4}  & circuit     & 67{,}029  & 1.381  & 7.3   & 5.5   & 124   & 2.9\\
\texttt{G2\_circuit} & circuit     & 150{,}102 & 4.030  & 1579  & 50.1  & 3011  & OOM\\
\bottomrule
\end{tabular}}
\end{table}

\section{Multicommodity congestion}\label{sec:mc}
\S\ref{sec:traffic} established the formulation and the $\Om$ solver for a single commodity. This
section extends both, end to end, to $K$ simultaneously routed commodities: a vector edge law that
couples the commodities through shared congestion (\S\ref{sec:mcform}), a solver that reuses
Algorithm~\ref{alg:fixed} and \emph{one} scalar inner setup (factor or hierarchy) for all $K$ commodities
(\S\ref{sec:mcalg}), validation against an exact-Newton baseline (full $Kn\times Kn$ block solve, \S\ref{sec:mcvalid}), and the same
full-corpus robustness sweep as \S\ref{sec:corpus} (\S\ref{sec:mccorpus}). The implementation is a
modular layer over the single-commodity solver, which remains intact.

\subsection{Formulation: a vector edge law}\label{sec:mcform}
$K$ commodities ship demands $\alpha d^1,\dots,\alpha d^K$ ($\mathbf 1^{\top}d^k=0$) over one
network; commodity $k$ has flow $f^k\in\R^m$ and potential $\phi^k\in\R^n$. Collect the per-edge
flows into the vector $\mathbf f_e=(f^1_e,\dots,f^K_e)\in\R^K$ and the per-edge potential gradients
into $\mathbf g_e=\big((B^{\top}\phi^1)_e,\dots,(B^{\top}\phi^K)_e\big)$. The modeling question is how the commodities couple through shared congestion. Three natural choices:
the signed total $\sum_k f^k_e$ (degenerate---collapses to single-commodity); the volume
$\sum_k|f^k_e|$ (classical but nonsmooth, forces a per-edge complementarity akin to the max-flow
cut); and the Euclidean magnitude $\|\mathbf f_e\|$ (smooth, strictly convex, reduces exactly to
$K=1$). We take the third for this feasibility study; the volume-coupled and directed models remain
future work (\S\ref{sec:multicommodity}). The program is \eqref{eq:primal} with the congestion cost charged to
the magnitude of the commodity-flow vector,
\begin{equation}\label{eq:mcprimal}
\min_{f^1,\dots,f^K}\ \sum_e\Phi_e\big(\|\mathbf f_e\|\big)\quad\text{s.t.}\quad
Bf^k=\alpha d^k,\quad k=1,\dots,K ,
\end{equation}
and stationarity gives, per edge, the \emph{vector law}
\begin{equation}\label{eq:mclaw}
\mathbf f_e=\rho_e\big(\|\mathbf g_e\|\big)\,\frac{\mathbf g_e}{\|\mathbf g_e\|}\,:
\end{equation}
the scalar law of \eqref{eq:law} applied to the magnitude of the gradient vector and acting along
it, so the $K$ gradient components share one congested conductance. At $K=1$, \eqref{eq:mclaw}
\emph{is} \eqref{eq:law}. The system $Bf^k=\alpha d^k$ with \eqref{eq:mclaw} is the stationarity of
the convex energy $E(\phi^1,\dots,\phi^K)=\sum_e\Psi_e(\|\mathbf g_e\|)-\alpha\sum_k(d^k)^{\top}\phi^k$,
strictly convex modulo the $K$ per-commodity constants, so the equilibrium is unique. Physically the coupling turns on exactly where congestion does. On real
capacity data at rush-hour load, commodity flows superpose to within solver tolerance ($\lesssim10^{-4}$ relative)---real networks place
each edge in one of BPR's two power-like regimes where path splits are load-independent---while the
shared load inflates marginal costs by up to $8.8\times$: the congestion externality is slowdown,
not rerouting. Synthetic instances with random capacities do reroute ($12$--$19\%$ flow deviation),
so the corpus sweep (\S\ref{sec:mccorpus}) exercises the solver where commodities genuinely
interact.

The Newton linearization of \eqref{eq:mclaw} is the $K$-commodity block Laplacian whose $(k,l)$
block is $B\,\mathrm{diag}\big(c_e\delta_{kl}+(\rho'_e-c_e)\,u^k_eu^l_e\big)B^{\top}$: per edge, the
$K\times K$ block
\begin{equation}\label{eq:mcjac}
\frac{\partial\mathbf f_e}{\partial\mathbf g_e}
= c_e\,I_K+(\rho'_e-c_e)\,\mathbf u_e\mathbf u_e^{\top},\qquad
c_e=\frac{\rho_e(s_e)}{s_e},\quad \mathbf u_e=\frac{\mathbf g_e}{s_e},\quad s_e=\|\mathbf g_e\| ,
\end{equation}
with eigenvalues $c_e$ (multiplicity $K-1$, across the flow direction) and $\rho'_e$ (along it):
symmetric positive definite for any monotone law, with anisotropy ratio
$c_e/\rho'_e=t_e'(F)F/t_e(F)\in[1,p]$ bounded by the BPR exponent.

\subsection{Algorithm: one scalar hierarchy carries all $K$ commodities}\label{sec:mcalg}
Algorithm~\ref{alg:fixed} carries over with three substitutions and no new machinery.
\emph{(i)~Stacked state.} The iterate is $\Phi=(\phi^1,\dots,\phi^K)$, the residual
$R=(r^1,\dots,r^K)$, $r^k=\alpha d^k-Bf^k$, and the stopping test is on the stacked norm.
\emph{(ii)~A scalar frozen operator.} Rather than freeze the $K\times K$ block linearization, NLF
freezes a scalar hierarchy $H$ on the single conductance
$w_e=\sqrt{c_e\rho'_e}$---the geometric mean of the block's two eigenvalues---and computes the
correction as $K$ independent inner solves, $\delta^k$ from $H$ against $r^k$. Thus one $\Om$ inner setup
serves all $K$ commodities and all chord steps.
Because BPR cost depends only on total edge flow, the $K\times K$ coupling per edge is rank-one, with
eigenvalues $c_e$ (multiplicity $K-1$) and $\rho'_e$; the $w_e$-preconditioned per-edge block then has
eigenvalues $\{\sqrt{c_e/\rho'_e},\sqrt{\rho'_e/c_e}\}$---symmetric about $1$ with ratio
$c_e/\rho'_e\le p$ independent of $K$, $m$, and the iterate. The per-edge conditioning is thus
uniformly controlled; the realized contraction of the assembled
iteration is governed by the graph-wide anisotropy and is measured at run time by the refresh monitor
below, not bounded a priori. \emph{(iii)~A self-calibrating refresh monitor.} The block-diagonal
frozen operator has an intrinsic rate floor set by the anisotropy that no rebuild can mitigate, so
a refresh monitor measures the contraction on the step following each rebuild and treats only
degradation beyond that baseline as staleness; otherwise a heavily congested instance would
rebuild every step to no effect.

The cost accounting is the point of the construction. A chord step costs \emph{one} vector-law
evaluation---$m$ scalar inversions, the same count as a single-commodity step, since only
$\|\mathbf g_e\|$ is inverted---plus $K$ cycles and $K$ residual products: $\mathcal{O}(Km)$, on top
of one shared $\Om$ setup. The entire $K$-dependence is $K$ right-hand sides through one hierarchy.

\subsection{Validation}\label{sec:mcvalid}
The frozen block Jacobian \eqref{eq:mcjac} matches the
finite-difference derivative of the residual; at $K=1$ both the multigrid and the direct paths
reproduce the single-commodity solver's equilibrium, and both agree with the exact-Newton baseline (full $Kn\times Kn$ block solve); identical commodities receive identical potentials and a negated
demand a negated potential (neither symmetry is imposed); and the joint-vs-alone coupling check
above. Table~\ref{tab:mcvalid} shows the reference comparison on real TNTP road networks at $K=4$
(four well-separated demand dipoles, benchmark load $\alpha=0.3\,\mathrm{median}(\{c_e\}_e)$ per commodity):
NLF reaches the exact-Newton equilibrium to $10^{-10}$ with a handful of chord steps and a
\emph{single} hierarchy setup---whereas the exact-Newton competitor must assemble and factor the full
$Kn\times Kn$ block linearization at \emph{every} step. This is the relevant high-accuracy baseline
for the coupled multicommodity program (general convex solvers such as Ipopt reduce to the same
the same $Kn\times Kn$ block linear algebra); NLF matches its equilibrium while replacing the $K$-fold-larger
block solve with one shared scalar $\Om$ hierarchy.

\begin{table}[t]\centering\small
\caption{Multicommodity validation on real road networks (TNTP), $K=4$. ``Exact Newton'' assembles
the full $Kn\times Kn$ block linearization and re-solves it at every damped step (the reference);
NLF is the block-diagonal chord with one frozen scalar hierarchy. Deviation is the relative
$\ell_2$ distance between the two equilibria.}\label{tab:mcvalid}
\begin{tabular}{lrrcccc}
\toprule
Network & $n$ & $m$ & Newton steps & NLF steps & Setups & Deviation\\
\midrule
SiouxFalls     & 24   & 38   & 8 & 8  & 1 & $5\times10^{-15}$\\
Anaheim        & 416  & 634  & 9 & 16 & 1 & $5\times10^{-10}$\\
Chicago-Sketch & 933  & 1475 & 9 & 19 & 1 & $3\times10^{-10}$\\
Winnipeg       & 1040 & 1595 & 6 & 19 & 1 & $7\times10^{-11}$\\
\bottomrule
\end{tabular}
\end{table}

\subsection{Applications across domains}\label{sec:mcapps}
We exercise the $K=4$ solver on the multicommodity applications the formulation targets: real TNTP
metropolitan road networks (topology, capacities, free-flow times, BPR parameters, and the four
heaviest origin--destination pairs all from the benchmark data, jointly scaled to a loaded rush hour,
peak $\|\mathbf f_e\|/c_e\approx1.07$), and corpus-protocol instances from other coupled-flow domains
where flows share one network---internet AS topology, a power transmission grid, national road
networks. Every instance converges to $10^{-9}$ with a handful of chord steps and $1$--$3$ hierarchy
setups, from $38$ edges to $1.5\times10^6$ (e.g.\ \texttt{roadNet-PA}, $1.5\times10^6$ edges, $26$
steps, $2$ setups).

\subsection{Corpus robustness and the cost of $K$}\label{sec:mccorpus}
\begin{figure}[t]
\centering
\includegraphics[width=\textwidth]{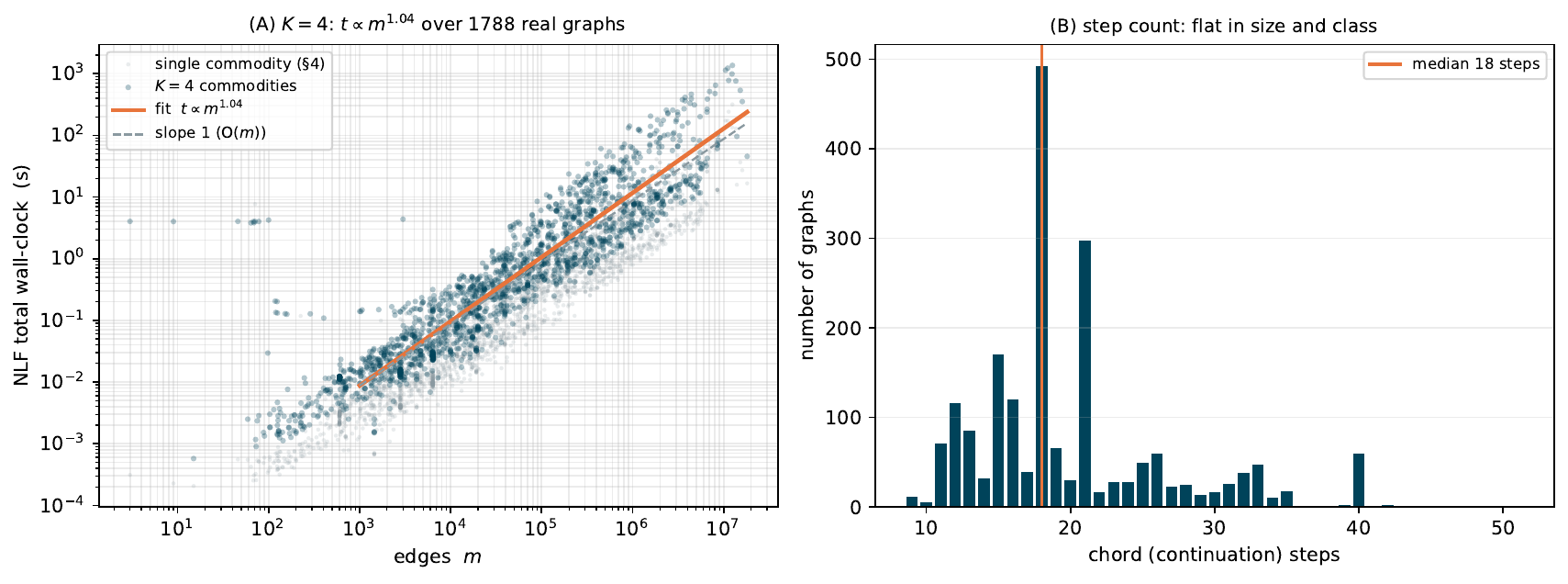}
\caption{Multicommodity NLF ($K=4$) over the full real-world SuiteSparse corpus, mirroring the
single-commodity sweep (grey, underneath). (A)~Total wall-clock vs.\ edge count: all $2{,}003$
graphs converge, $t\propto m^{1.04}$ over four decades. (B)~Chord-step count: flat in size and
class, median $18$---unchanged from the single-commodity solver.}\label{fig:mccorpus}
\end{figure}

The decisive test is the same full-corpus sweep as \S\ref{sec:corpus}, repeated at $K=4$: the
identical $2{,}003$ real-world SuiteSparse graphs and instance protocol, with four well-separated
demand dipoles per graph, each at the benchmark load. Every one of the $2{,}003$ graphs
converges to $10^{-9}$ relative residual (Fig.~\ref{fig:mccorpus}), to $m=1.8\times10^{7}$
edges, at a flat median of $18$ chord steps (range $9$--$51$) and a median of $1$ hierarchy
setup (maximum $5$; more than two on only $120$ of $2{,}003$, triggered by the refresh conditions of \S\ref{sec:inner}): the shared scalar hierarchy carries
the coupled four-commodity system essentially as cheaply as it carries one. The fitted total cost
is $t\propto m^{1.04}$ over four decades, median $9.1~\mu$s per edge. Joined per graph against the
single-commodity sweep, the $K=4$ solve costs a median $3.8\times$ the single-commodity solve
(interquartile range $[3.0,5.6]$)---the $\approx\!K$ per-step work with the setup amortized, as
the cost accounting of \S\ref{sec:mcalg} predicts.

The cost of $K$ is measured directly by sweeping $K\in\{1,2,4,8\}$ with nested demand sets on FEM and
road graphs to $1.5\times10^6$ edges (\texttt{olesnik0}, \texttt{srb1}, \texttt{roadNet-PA}): the
chord-step count and the hierarchy-setup count are roughly constant, and wall-clock grows close to linearly in $K$ ($K{=}8$ at $4.2$--$8.5\times$ the single-commodity cost), i.e., $\mathcal{O}(Km)$ in total.

\section{Extensions and Future work}\label{sec:extensions}

\subsection{Directed multicommodity assignment}\label{sec:multicommodity}
\S\ref{sec:mc} routes $K$ coupled commodities under a smooth magnitude coupling with signed flows;
real traffic assignment is \emph{directed} (nonnegative arc flows, costs coupling through the directed
link flow), validated against published origin--destination equilibria (Sioux Falls, Anaheim, \dots)~\cite{tntp,bargera2002}.
Two steps separate \S\ref{sec:mc} from that target: the volume coupling $\sum_k|f^k_e|$ forces a
per-edge complementarity (only commodities whose potential drop matches the common congested cost may
flow)---a free boundary of the same character as the max-flow cut, attacked by the
continuation-with-safeguard machinery of \S\ref{sec:continuation}; and the scale of $K$ (thousands of
demand columns), made conceivable by the $\mathcal{O}(Km)$ per-step cost with one shared hierarchy
(\S\ref{sec:mcalg}) plus per-origin aggregation borrowed from the bush-method literature~\cite{bargera2002,bargera2010}, making the $K$ scaling conceivable. Benchmarking against Frank--Wolfe and
Algorithm~B on the published equilibria would be required.

\paragraph{The payoff: design optimization} The forward equilibrium is the means, not the end. Its
motivating problems---congestion pricing, origin--destination matrix calibration,
capacity network design---are bilevel programs optimized over a design $x$ (tolls, demand,
capacities) subject to the equilibrium $f(x)$, where the load-bearing object is the sensitivity
$\partial f/\partial x$. Because $f$ is the stationarity of \eqref{eq:master}, that sensitivity is a
single \emph{adjoint} solve with the \emph{same} symmetric Laplacian Jacobian: the gradient of an outer
objective with respect to all $k$ design variables costs one additional $\Om$ solve,
independent of $k$. Equilibrium sensitivity is itself classical---the adjoint/implicit-function
calculus of Tobin and Friesz~\cite{tobin-friesz} and the bilevel network-design literature it seeds
compute exactly this derivative---so the contribution is not the adjoint but that here it \emph{is} the
same $\Om$ Laplacian solve: available natively, and near-linear rather than a direct KKT factorization
on poorly-separable designs. A directed single-commodity prototype we developed validates this on
optimal tolling: the adjoint gradient matches finite differences. The undirected of this paper is thus a stepping stone: the same
$\Om$ Laplacian machinery that \emph{solves} \eqref{eq:master} also \emph{differentiates} it, which is
exactly what bilevel network design needs.

\subsection{One-shot FAS}\label{sec:fasfail}
The methodologically pure alternative to Algorithm~\ref{alg:fixed} is a Full Approximation Scheme (FAS) multilevel cycle (local linearization)~\cite{brandt-fas,brandt-guide}. In initial experiments, the LAMG+ caliber-1 interpolation misrepresented the flux feeding the coarse $\tau$ correction: caliber-1 (piecewise-constant) prolongation implicitly assumes the fine-level conductances $\rho'_e$ are near-uniform within an aggregate, but near saturation $\rho'$ varies by orders of magnitude across a single aggregate, so the Galerkin coarse operator $P^{\top}JP$ no longer represents the fine-level flux balance and the nonlinear coarse correction is inconsistent. Building a variationally faithful transfer for the saturating law---e.g.\ a flux-based interpolation~\cite[\S4.6]{brandt-guide}---remains an open question.

\subsection{Full multigrid: $\alpha$ update at the coarsest level}\label{sec:fmg}
Brandt's nested-iteration ideal~\cite[\S5.6]{brandt-guide} would update the global scalar $\alpha$ at
the coarsest level rather than the finest, prolongating $(\phi,\alpha)$ upward with per-level
arclength finishing. A first version was $1.6$--$2.1\times$ faster on small instances but $2.2\times$
slower on a larger grid (per-level finishing dominates)---a non-uniform gain---so we keep the simple
$\Om$ finest-level architecture and leave a tuned FMG variant to future work.

\subsection{Optimal transport as a further instance}\label{sec:ot}
The same edge-separable energy expresses optimal transport in Beckmann (flow) form~\cite{beckmann}: minimizing a
monotone edge-separable cost subject to $Bf=\mu-\nu$ is exactly \eqref{eq:master}.
The $W_p$ Wasserstein distance~\cite{peyre-cuturi} between measures $\mu$ and $\nu$ on a graph
corresponds to the $\ell^p$ edge cost; the quadratic case ($p=2$, $H^{-1}$ cost) reduces to a single $\Om$ Laplacian solve, while the $W_1$ limit ($p\to1$) gives the non-smooth saturating
member, a natural fold instance. The quadratically-regularized form, whose dual Hessian is an
active-subgraph Laplacian solved in the literature by direct factorization~\cite{essid-solomon}, is
where a near-linear inner could help; whether it is competitive is an open question, since the
$W_1$/Beckmann limit is a min-cost flow with strong combinatorial and direct solvers.

\subsection{DC optimal power flow}\label{sec:dcopt}

The DC approximation to optimal power flow (DC-OPF) minimizes generation cost
subject to power balance and thermal line limits:
\begin{equation}\label{eq:dcopf}
\begin{aligned}
  &\min_{\theta,\,p^g}\;\sum_i C_i(p^g_i)
  \quad\text{s.t.}\quad
  B\,\mathrm{diag}(b)\,B^{\top}\theta = p^g - p^d,\\
  &\qquad\quad
  |b_e(B^{\top}\theta)_e|\le c_e,\quad
  p^{\min}\le p^g\le p^{\max},
\end{aligned}
\end{equation}
where $b_e$ is the line susceptance, $p^d$ fixed demand, and $C_i$ a convex
(typically quadratic or piecewise-linear) generation cost.
The power-balance constraint is a susceptance-weighted graph Laplacian; the
line limits are the box constraints of the saturating law~\eqref{eq:rho},
so \eqref{eq:dcopf} is an instance of \eqref{eq:master} extended by a nodal
cost and generator box constraints.

We merely record DC-OPF as an instance of the framework but do not solve the constrained problem here: the
generator box constraints and nodal cost $\sum_i C_i(p^g_i)$ require an FAS extension~\cite{brandt-fas}
that carries $p^g$ and its bounds into the coarse levels (\S\ref{sec:fasfail}), which we have not built.
As a back-of-envelope projection only---not a measured result---current IPM solvers (MATPOWER,
PowerModels.jl) run $20$--$50$ Newton steps, each a fresh sparse-direct KKT factorization; on
near-planar grids ($118$ to $30{,}000$ buses) we timed that factorization at $2$--$6\times$ cheaper than one
NLF cycle, so whether a continuation-Newton loop reusing one frozen hierarchy would win on total cost
turns on the iteration-count trade and is an open question.

\section{Conclusion}\label{sec:conclusion}
Placing the nonlinearity in the \emph{edge law} collapses a family of convex, edge-separable
network-flow equilibria---electrical, congestion, minimum-delay, maximum-flow---onto one nonlinear
Laplacian $B\,\rho(B^{\top}\phi)=\alpha d$, solved by a damped chord-Newton on a frozen linearization
cycled against the nonlinear residual. Any
near-linear Laplacian solver is shown to empirically translate into the same complexity on the whole nonlinear class.

\paragraph{Reproducibility} Every table is produced by released scripts, on public benchmarks---the
SuiteSparse collection (each instance reduced to its largest connected component, \S\ref{sec:corpus}),
the Transportation Networks road data~\cite{tntp}, and MATPOWER grids~\cite{matpower}---and on
seeded synthetic generators (Erd\H{o}s--R\'enyi topologies for the scaling study) where noted. The interior-point
comparator is Ipopt~3.14.19 (MUMPS~5.9.0, $120$\,s CPU budget); the BPR law uses the strictly convex
regularization of \S\ref{sec:bpr}; capacity/free-flow distributions and demand seeds are fixed and
recorded. The inner solvers are public~\cite{approxchol,lamgplus}; the NLF solver, instance-generation
scripts, and a reproducibility notebook that re-runs the timings and checks them against these tables
are at \url{https://github.com/orenlivne/nlf}. All timings were collected on a MacBook Pro with Apple M5 Pro (18-core: 6 efficiency + 12 performance) and 48\,GB unified memory, running a single Julia process with no explicit parallelism.

\appendix
\section{Why each application is an instance of the framework}\label{app:instances}
Per application: equilibrium principle $\Leftrightarrow$ convex program \eqref{eq:primal}
$\Leftrightarrow$ KKT $\Leftrightarrow$ $B\rho(B^{\top}\phi)=\alpha d$.
Each instance differs only in the edge law; directed complementarity is absorbed into that law or
relaxed by the signed formulation.

\paragraph{The chain in general} Stationarity of \eqref{eq:primal} in $f_e$ gives the nonlinear
Ohm law $t_e(f_e)=(B^{\top}\phi)_e$; inversion gives $f_e=\rho_e((B^{\top}\phi)_e)$; and
Kirchhoff ($Bf=\alpha d$) yields \eqref{eq:master}. Summing along any source--sink path telescopes:
every path has marginal cost $\phi_s-\phi_t$, the scalar equilibrium cost. The signed formulation
has no inequality constraints---a smooth equation---which is why Newton with a near-linear Laplacian
solve applies with no active sets.

\paragraph{Electrical flow (linear law)} Ohm's law is linear, $f_e=g_e/r_e$; the energy is the
dissipated power and \eqref{eq:primal} is Thomson's principle (current minimizes
dissipation)~\cite{doyle-snell}. Equation \eqref{eq:master} is the ordinary weighted Laplacian
system---the framework's fixed point under linearization, and the object the inner solver
\cite{lamgplus} is built for.

\paragraph{Congestion (traffic) equilibrium} Wardrop's principle---per origin--destination pair,
path flows $h_p\ge0$ with $\sum_p h_p=d$ obeying $T_p-\pi\ge0$ and $h_p(T_p-\pi)=0$---is, by
Beckmann's theorem~\cite{beckmann}, the KKT system of $\min_h\sum_e\int_0^{f_e}t_e$: differentiating
gives $\partial F/\partial h_p=\sum_{e\in p}t_e=T_p$, and $\pi$ (the multiplier of $\sum_p h_p=d$) is
both the marginal cost of demand and the equilibrium travel time. In link variables the demand
multiplier unfolds into the node field $\phi$ ($\pi=\phi_s-\phi_t$) and the program becomes
\eqref{eq:primal} with $\Phi_e'=t_e$, the BPR latency of \S\ref{sec:traffic}. The signed relaxation
keeps the first Wardrop clause exactly (all routes at potential-drop cost $\pi$) and discards
$h_p\ge0$: an unused route may carry counterflow---the undirected scope, with per-commodity
nonnegativity the free-boundary extension of \S\ref{sec:extensions}.

\paragraph{Minimum-delay routing} M/M/1 queueing with Little's law gives the per-packet link
delay $t_e(f)=1/(c_e-f)$, and Kleinrock's independence approximation makes total delay
edge-separable~\cite{bertsekas-gallager}. Taking the delay as marginal cost,
$\Phi_e=\log\!\big(c_e/(c_e-f)\big)$, gives \eqref{eq:primal} directly; the inverse law
$\rho_e(g)=c_e-1/g$ (odd-extended) saturates at capacity with conductance
$\rho'_e=(c_e-f)^2\to0$, placing the problem in the fold class of Table~\ref{tab:instances}: the
capacity constraint $f<c_e$ is never imposed---the law cannot produce it---and the feasibility
limit, where the filling edges form a cut in the spare capacities, is the fold handled by
\S\ref{sec:continuation}.

\paragraph{Maximum flow} The LP $\max\,\alpha$ s.t.\ $Bf=\alpha d$, $|f_e|\le c_e$ has the box
constraint as its only nonsmoothness; its complementarity (an edge is either slack with zero
reduced cost or saturated on the cut) is absorbed into the smooth saturating law \eqref{eq:rho},
whose energy $\Psi_e=c_e^2\log\cosh(g/c_e)$ is quadratic for slack edges and exactly linear,
$c_e|g|$, for saturated ones. The active set is not tracked: it \emph{emerges} as the
zero-conductance set $w_e\to0$, the forming min cut, and the LP's degeneracy reappears as the one
scalar fold at $\alpha=F^{*}$ traversed by the continuation of \S\ref{sec:continuation}. At the
fold the saturated set is the exact min cut and $\phi$ steps across it (\S\ref{sec:maxflow}); no
$(1+\epsilon)$ relaxation is involved.

\end{document}